\documentclass[a4paper,11pt] {amsart}
\usepackage{amsthm,amssymb,latexsym,mathrsfs}
\input amssym.def

\author{Beno\^it F. Sehba}

\title[Maximal functions and Fractional Bergman operators ]{Weighted boundedness of maximal functions and fractional Bergman operators}
\newtheorem{theorem}{T{\hskip 0pt\footnotesize\bf HEOREM}}[section]
\newtheorem{lemma}[theorem]{L{\hskip 0pt\footnotesize\bf EMMA}}
\newtheorem{proposition}[theorem]{P{\hskip 0pt\footnotesize\bf ROPOSITION}}
\newtheorem{definition}[theorem]{D{\hskip 0pt\footnotesize\bf EFINITION}}
\newtheorem{corollary}[theorem]{C{\hskip 0pt\footnotesize\bf OROLLARY}}

\newcommand{\bprop} {\begin{proposition}}
\newcommand{\eprop} {\end{proposition}}
\newcommand{\btheo} {\begin{theorem}}
\newcommand{\etheo} {\end{theorem}}
\newcommand{\blem} {\begin{lemma}}
\newcommand{\elem} {\end{lemma}}
\newcommand{\bcor} {\begin{corollary}}
\newcommand{\ecor} {\end{corollary}}

\newcommand{\Be}{\begin{equation}}
\newcommand{\Ee}{\end{equation}}
\newcommand{\Bea}{\begin{eqnarray}}
\newcommand{\Eea}{\end{eqnarray}}
\newcommand{\Bes}{\begin{equation*}}
\newcommand{\Ees}{\end{equation*}}
\newcommand{\Beas}{\begin{eqnarray*}}
\newcommand{\Eeas}{\end{eqnarray*}}
\newcommand{\Ba}{\begin{array}}
\newcommand{\Ea}{\end{array}}

\scrollmode
\begin{document}
\address{Beno\^it F. Sehba, Department of Mathematics, University of Ghana, Legon, P. O. Box LG 62, Legon, Accra, Ghana}
\email{bfsehba@ug.edu.gh}


\maketitle

\begin{abstract}
The aim of this paper is to study two-weight norm inequalities for
fractional maximal functions and fractional Bergman operator defined on the upper-half space. Namely, we characterize those pairs of weights for which these maximal operators satisfy strong and weak
type inequalities. Our characterizations are in terms of Sawyer and B\'ekoll\'e-Bonami type conditions. We  also obtain a “$\Phi$-bump” characterization for these maximal functions, where $\Phi$ is a Orlicz function.  As a consequence, we obtain two-weight norm inequalities for fractional Bergman operators. Finally, we provide some sharp weighted inequalities for the fractional maximal functions.
\end{abstract}


\section{Introduction}
Let $\mathcal H$ be the upper-half plane, that is the set $\{z=x+iy\in \mathbb {C}: x\in \mathbb R,\,\,\, \textrm{and}\,\,\,y>0\}$. Given $\omega$ a nonnegative locally integrable function on $\mathcal H$ (i.e. a weight), $\alpha>-1$, and $1\le p<\infty$, we denote by $L^p(\mathcal H, \omega dV_\alpha)$, the set of functions $f$ defined on $\mathcal H$ such that
$$||f||_{p,\omega, \alpha }^p:=\int_{\mathcal{H}}|f(z)|^p\omega(z)dV_\alpha(z)<\infty$$
with $dV_\alpha(x+iy)=y^\alpha dxdy$. We write $L^p(\mathcal H, dV_\alpha)$ when $\omega(z)=1$ for any $z\in \mathcal{H}$ and $\|\cdot\|_{p,\alpha}$ for the corresponding norm.
\vskip .1cm
For $\alpha>-1$ and $0\le \gamma<2+\alpha$, the positive fractional Bergman operator $T_{\alpha,\gamma}$ is defined by
\Be\label{eq:fracBergdef}
T_{\alpha,\gamma}f(z):=\int_{\mathcal{H}}\frac{f(w)}{|z-\overline{w}|^{2+\alpha-\gamma}}dV_\alpha(w).
\Ee
For $\gamma=0$, the operator $T_\alpha:=T_{\alpha, 0}$ corresponds to the positive Bergman projection.
\vskip .3cm
For $I\subset \mathbb{R}$ an interval, we put $Q_I:=\{z=x+iy\in \mathbb C:x\in I\,\,\,\textrm{and}\,\,\,0<y<|I|\}$. The fractional maximal function $\mathcal{M}_{\alpha,\gamma} $ is the function defined for any $f\in L_{loc}^1(\mathcal H)$ by
$$\mathcal{M}_{\alpha,\gamma} f(z)=\sup_{I\in \mathcal {D}^\beta}\frac{\chi_{Q_I}(z)}{|I|^{2+\alpha-\gamma}}\int_{Q_I}|f(w)|dV_\alpha(w).$$
When $\gamma=0$, $\mathcal{M}_{\alpha}:=\mathcal{M}_{\alpha,0}$ is just the Hardy-Littlewood maximal function.

The operators $T_{\alpha,\gamma}$ and $\mathcal{M}_{\alpha,\gamma}$ appear naturally in the problem of off-diagonal weighted inequalities for the Bergman operator (see \cite{Sehba2}). Obviously, $\mathcal{M}_{\alpha,\gamma}$ is pointwise dominated by $T_{\alpha,\gamma}$ and it is easy to see that given two weights $\sigma$ and $\omega$ on $\mathcal{H}$, for $1\le p\le q<\infty$, the boundedness of $T_{\alpha,\gamma}$ from $L^p(\mathcal H, \sigma dV_\alpha)$ to $L^q(\mathcal H, \omega dV_\alpha)$ implies the boundedness of $P_\alpha^+=T_{\alpha}$ from  $L^p(\mathcal H, \sigma dV_\alpha)$ to $L^q(\mathcal H, \omega dV_\beta)$ where $\beta=\alpha+q\gamma $.
\vskip .3cm
We are interested in this work to the pairs of measure $\mu$ and weight $\omega$ such that the operator $\mathcal{M}_{\alpha,\gamma}$ satisfies strong and weak type inequalities. More precisely, given $1\le p\le q<\infty$, we provide some characterizations of positive measures $\mu$ on $\mathcal H$ and weight $\omega$ such that the following strong inequality holds
\Be\label{eq:strongineq1}
\int_{\mathcal H}\left(\mathcal{M}_{\alpha,\gamma} f(z)|\right)^qd\mu(z)\le C\left(\int_{\mathcal H}|f(z)|^p\omega(z)dV_\alpha(z)\right)^{\frac{q}{p}}
\Ee

\vskip .1cm
We also characterize those positive measures $\mu$ on $\mathcal H$ such that
\Be\label{eq:weakineq1}
\mu\left(\{z\in \mathcal H: \mathcal{M}_{\alpha,\gamma}f(z)>\lambda\}\right)\le \frac{C}{\lambda^q}\left(\int_{\mathcal H}|f(z)|^p\omega(z)dV(z)\right)^{\frac{q}{p}}
\Ee

In the case of strong inequalities, our characterizations are given in terms of Sawyer type conditions, B\'ekoll\`e-Bonami type conditions. Sufficient conditions are also obtained by adding some $\Phi$-bump conditions on the weight, where $\Phi$ is an appropriate Young function. The latter allows us to obtain two-weight norm inequalities for the fractional Bergman operators. Finally, we provide some weighted norm estimates for the above fractional maximal function, we prove that some of these estimates are sharp.
\vskip .3cm
\section{Statement of the results}
Let $\alpha>-1$ and let $\omega$ be a weight. For any subset $E$ of $\mathcal{H}$, we use the notation
$|E|_{\omega,\alpha}:=\int_E \omega(z)dV_\alpha(z)$. When $\omega=1$, we simply write $|E|_\alpha$.
Let us recall that for any interval $I\subset \mathbb{R}$, its associated Carleson square $Q_I$ is the set $$Q_I:=\{z=x+iy\in \mathbb C:x\in I\,\,\,\textrm{and}\,\,\,0<y<|I|\}.$$

Let $\alpha>-1$, and $1<p<\infty$. Given a weight $\omega$, we say $\omega$ is in the B\'ekoll\`e-Bonami class $\mathcal{B}_{p,\alpha}$,  if the quantity
$$[\omega]_{\mathcal{B}_{p,\alpha}}:=\sup_{I\subset \mathbb R }\left(\frac{1}{|Q_I|_\alpha}\int_{Q_I}\omega(z)dV_\alpha(z)\right)\left(\frac{1}{|Q_I|_\alpha}\int_{Q_I}\omega(z)^{1-p'}dV_\alpha(z)\right)^{p-1}$$
is finite. This is the exact range of weights $\omega$ for which the orthogonal projection $P_\alpha$ from $L^2(\mathcal H, dV_\alpha(z))$ onto its closed subspace consisting of analytic functions is bounded on $L^p(\mathcal H, \omega dV_\alpha)$ (see \cite{Bek,BB,PR}). For $p=\infty$, we say $\omega\in \mathcal{B}_{\infty,\alpha}$, if 
$$[\omega]_{\mathcal{B}_{\infty,\alpha}}:=\sup_{I\subset \mathbb R }\frac{1}{|Q_I|_{\omega,\alpha}}\int_{Q_I}\mathcal{M}_\alpha(\omega\chi_{Q_I})dV_\alpha(z)<\infty.$$

\subsection{Some weak inequalities}
Our first result gives some elementary (unweighted) inequalities for the fractional maximal function.
\btheo\label{thm:main1}
Let $\alpha>-1$, and $0\le \gamma<2+\alpha$. Then the following hold.
\begin{itemize}
\item[(a)] For any $\lambda>0$, there exists a positive constant $C$ such that
\Be\label{eq:main11}|\{z\in \mathcal{H}:\mathcal{M}_{\alpha,\gamma}f(z)>\lambda\}|_\alpha\le C\left(\frac{1}{\lambda}\int_{\mathcal{H}}|f(z)|dV_\alpha(z)\right)^{\frac{2+\alpha}{2+\alpha-\gamma}}.
\Ee
\item[(b)] $$\sup_{z\in \mathcal{H}}\mathcal{M}_{\alpha,\gamma}f(z)\le \|f\|_{\frac{2+\alpha}{\gamma},\alpha}$$
where for $\gamma=0$, $ \|f\|_{\frac{2+\alpha}{\gamma},\alpha}$ is understood as $ \|f\|_{\infty}$.
\end{itemize}
\etheo
It follows from the above result and Marcinkiewicz interpolation theorem that the following holds.
\bcor\label{cor:main1}
Let $\alpha>-1$, and $0\le \gamma<2+\alpha$. Then $\mathcal{H}:\mathcal{M}_{\alpha,\gamma}$ is bounded from $L^p(\mathcal{H}, dV_\alpha)$ to $L^q(\mathcal{H}, dV_\alpha)$, for $1<p<\frac{2+\alpha}{\gamma}$ and $\frac{1}{p}-\frac{1}{q}=\frac{\gamma}{2+\alpha}$.
\ecor
Our next result provides weak-type estimates.
\btheo\label{thm:main2}
Let $1\le p\le q<\infty$, $\alpha>-1$ and $0\le \gamma<2+\alpha$. Let $\omega$ be a weight on $\mathcal H$. Then the following assertions are equivalent.
\begin{itemize}
\item[(a)] There is a constant $C_1>0$ such that for any $f\in L^p(\mathcal H, \omega dV_\alpha)$, and any $\lambda>0$,
\Be\label{eq:main21}
\mu(\{z\in \mathcal H: \mathcal{M}_{\alpha,\gamma}f(z)>\lambda\})\le \frac{C_1}{\lambda^q}\left(\int_{\mathcal H}|f(z)|^p\omega(z)dV_\alpha(z)\right)^{q/p}
\Ee
\item[(b)] There is a constant $C_2>0$ such that for any interval $I\subset \mathbb R$,
\Be\label{eq:main22}
|Q_I|_\alpha^{q(\frac{\gamma}{2+\alpha}-\frac{1}{p})}\left(\frac{1}{|Q_I|_\alpha}\int_{Q_I}\omega^{1-p'}(z)dV_\alpha(z)\right)^{q/p'}\mu(Q_I)\le C_2
\Ee
where $\left(\frac{1}{|Q_I|_\alpha}\int_{Q_I}\omega^{1-p'}(z)dV_\alpha(z)\right)^{1/p'}$ is understood as $\left(\inf_{Q_I}\omega\right)^{-1}$ when $p=1$.
\item[(c)] There exists a constant $C_3>0$ such that for any locally integrable function $f$ and  any interval $I\subset \mathbb R$,
\Be\label{eq:main23}
\left(\frac{1}{|Q_I|_\alpha^{1-\frac{\gamma}{2+\alpha}}}\int_{Q_I}|f(z)|dV_\alpha(z)\right)^q\mu(Q_I)\le C_3\left(\int_{Q_I}|f(z)|^p\omega(z) dV_\alpha(z)\right)^{q/p}.
\Ee
\end{itemize}
\etheo
\subsection{Strong inequalities}
We also observe the following Sawyer-type characterization.
\btheo\label{thm:main3}
Let $1<p\le q<\infty$, $\alpha>-1$, and $0\le \gamma<2+\alpha$. Let $\mu$ be a positive measure and $\sigma$ a weight on $\mathcal{H}$. Then the following are equivalent.
\begin{itemize}
\item[(a)] There exists a positive constant $C_1$ such that
\Be\label{eq:main31}\left(\int_{\mathcal{H}}\left(\mathcal{M}_{\alpha,\gamma}(\sigma f)(z)\right)^qd\mu(z)\right)^{1/q}\le C_1\left(\int_{\mathcal{H}}|f(z)|^p\sigma(z)dV_\alpha(z)\right)^{1/p}.
\Ee
\item[(b)] There is a positive constant $C_2$ such that for any interval $I\subset \mathbb{R}$,
\Be\label{eq:main32}\left(\int_{\mathcal{H}}\left(\mathcal{M}_{\alpha,\gamma}(\chi_{Q_I}\sigma )(z)\right)^qd\mu(z)\right)^{1/q}\le C_2|Q_I|_{\sigma,\alpha}^{1/p}.
\Ee
\end{itemize}
Moreover,
$$\|\mathcal{M}_{\alpha,\gamma} (\sigma f)\|_{L^q{d\mu}}\simeq \sup_{I\subset \mathbb{R}}\left(\frac{\left(\int_{\mathcal{H}}\left(\mathcal{M}_{\alpha,\gamma}(\chi_{Q_I}\sigma )(z)\right)^qd\mu(z)\right)^{1/q}}{|Q_I|_{\sigma,\alpha}^{1/p}}\right).$$
\etheo
We have the following result for the strong inequality.
\btheo\label{thm:main4}
Let $\alpha>-1$, $1<p\le q<\infty$, and and $0\le \gamma<2+\alpha$. Let $\mu$ be a positive measure and $\sigma$ a weight on $\mathcal{H}$. Assume that $\sigma\in \mathcal{B}_{\infty,\alpha}$. Then the following assertions are equivalent.
 \begin{itemize}
\item[(i)] There exists a constant $C_1>0$ such that for any $f\in L^p(\mathcal H, \sigma dV_\alpha)$,
\Be\label{eq:main41}
\left(\int_{\mathcal H}\left(\mathcal{M}_{\alpha,\gamma} \sigma(z)f(z)\right)^qd\mu(z)\right)^{1/q}\le C_1\|f\|_{p,\sigma,\alpha}.
\Ee

\item[(ii)] There is a constant $C_2$ such that 
\Be\label{eq:main42} [\sigma,\mu]_{p,q.\alpha,\gamma}:=\sup_{I\subset \mathbb{R}}|Q_I|_{\alpha}^{-q(1-\frac{\gamma}{2+\alpha})}\mu(Q_I)|Q_I|_{\sigma,\alpha}^{\frac{q}{p'}}\leq C_2.\Ee
\end{itemize}
Moreover, $$\|\mathcal{M}_{\alpha,\gamma} (\sigma f)\|_{L^q{d\mu}}\le ([\sigma,\mu]_{p,q.\alpha,\gamma})^{1/q}[\sigma]_{\mathcal{B}_{\infty,\alpha}}^{1/p}\|f\|_{p,\sigma,\alpha}.$$
\etheo
\subsection{Bump-condtion for the fractional operators}
Recall that a function from $[0,\infty)$ to itself is a Young function if it is continuous, convex and increasing, and satisfies $\Phi(0)=0$ and $\Phi(t)\rightarrow \infty$ as $t\rightarrow \infty$. 
\vskip .3cm
Given a Young function $\Phi$, we say it satisfies the $\Delta_2$ (or doubling) condition, if there exists a constant $K>1$ such that, for any $t\ge 0$,
\begin{equation}\label{eq:delta2condition}
 \Phi(2t)\le K\Phi(t).
 \end{equation}

Let $\Phi$ be a Young function, and $\alpha>-1$. For any interval $I\subset \mathbb{R}$, define $L^\Phi(Q_I, |Q_I|_{\alpha}^{-1}dV_\alpha)$ to be the space of all functions $f$ such that 
$$\frac{1}{|Q_I|_{\alpha}}\int_{Q_I}\Phi\left(|f(z)|\right)dV_\alpha(z)<\infty.$$
We define on $L^\Phi(Q_I, |Q_I|_{\alpha}^{-1}dV_\alpha)$ the following Luxembourg norm $$\|f\|_{Q_I,\Phi,\alpha}:=\inf\{\lambda>0: \frac{1}{|Q_I|_{\alpha}}\int_{Q_I}\Phi\left(\frac{|f(z)|}{\lambda}\right)dV_\alpha(z)\le 1\}.$$
 When $\Phi(t)=t^p$, $1\le p<\infty$, $L^\Phi(Q_I, |Q_I|_{\alpha}^{-1} dV_\alpha)$ is just $L^p(Q_I, |Q_I|_{\alpha}^{-1}dV_\alpha)$ in which case $\|f\|_{Q_I,\Phi,\alpha}$ is just replaced by $\left(\frac{1}{|Q_I|_{\alpha}}\int_{Q_I}|f(z)|^pdV_\alpha(z)\right)^{1/p}$. Then the maximal function $\mathcal{M}_{\Phi,\alpha}$ is defined as $$\mathcal{M}_{\Phi,\alpha}f(z):=\sup_{I\subset \mathbb{R},z\in Q_I}\|f\|_{Q_I,\Phi,\alpha}.$$

We recall that the complementary function $\Psi$ of the Young function $\Phi$, is the function defined from $\mathbb R_+$ onto itself by
\begin{equation}\label{complementarydefinition}
\Psi(s)=\sup_{t\in\mathbb R_+}\{ts - \Phi(t)\}.
\end{equation}
Let $1<p<\infty$. We say a Young function $\Phi$ belongs to the class $B_p$, if it satisfies the $\Delta_2$ condition and there is a positive constant $c$ such that 
\Be\label{eq:BpDini}
\int_c^\infty \frac{\Phi(t)}{t^p}\frac{dt}{t}<\infty.
\Ee
The following result provides a sufficient condition for the off-diagonal boundedness of the maximal function $\mathcal{M}_{\alpha,\gamma}$.
\btheo\label{thm:main5}
Let $\alpha>-1$, $1<p\le q<\infty$ and $0\le \gamma<2+\alpha$. Let $\Phi\in B_p$ and denote by $\Psi$ its complementary function. Assume that $\omega$ is a weight and $\mu$ is a positive Borel measure on $\mathcal H$ such that there is a positive constant $C$ for which for any interval $I\subset \mathbb{R}$,
\Be\label{eq:main51}
|Q_I|_\alpha^{q(\frac{\gamma}{2+\alpha}-\frac{1}{p})}\|\omega^{-1}\|_{Q_I,\Psi,\alpha}^q\mu(Q_I)\le C.
\Ee

 Then there is a positive constant $K$ such that for any $f\in L^p(\mathcal H, \omega dV_\alpha)$, 
 \Be\label{eq:main52}
 \left(\int_{\mathcal H}\left(\mathcal{M}_{\alpha,\gamma} f(z)\right)^qd\mu(z)\right)^{1/q}\le K\|f\omega\|_{p,\alpha}.
 \Ee
\etheo
Conditions of type (\ref{eq:main51}) are known as bump-conditions. The above result is used to prove the following sufficient condition for the boundedness of the fractional Bergman operator.
 \btheo\label{thm:main6}
Suppose $\alpha>-1$, $0\le \gamma<2+\alpha$, and $1<p\le q<\infty$. Let $\Phi$ and $\Psi$ be two Young functions whose complementary functions are respectively in $B_p$ and $B_{q'}$. Assume that $\omega$ and $\sigma$ are weights on $\mathcal H$ such that there is positive constant $C$ for which for any interval $I\subset \mathbb{R}$,
\Be\label{eq:main61}
|Q_I|_\alpha^{\frac{\gamma}{2+\alpha}+\frac{1}{q}-\frac{1}{p}}\|\omega\|_{Q_I,\Psi,\alpha}\|\sigma^{-1}\|_{Q_I,\Phi,\alpha}\le C.
\Ee

 Then there is a positive constant $K$ such that for any $f\in L^p(\mathcal H, \sigma dV_\alpha)$, 
 \Be\label{eq:main62}
 \left(\int_{\mathcal H}\left(\omega(z)T_{\alpha,\gamma} f(z)\right)^qdV_\alpha(z)\right)^{1/q}\le K\|f\sigma\|_{p,\alpha}.
 \Ee
\etheo
\subsection{Some weighted norm inequalities for $\mathcal{M}_{\alpha,\gamma}$}
Let $1<p,q<\infty$, $\alpha>-1$, and $0\le \gamma<2+\alpha$. We introduce the classes $\mathcal{A}_{p,q,\alpha}$, $\mathcal{C}_{p,q,\alpha}$, and $\mathcal{S}_{p,q,\alpha}$ of pairs of weights. We say the pair of weights $(\sigma,\omega)$ belongs to $\mathcal{A}_{p,q,\alpha}$, if 
\Be\label{eq:Apqclass}
[\sigma,\omega]_{\mathcal{A}_{p,q,\alpha}}:=\frac{|Q_I|_{\omega,\alpha}^{p/q}|Q_I|_{\sigma,\alpha}^{p/p'}}{|Q_I|_\alpha^{p(1-\frac{\gamma}{2+\alpha})}}<\infty.
\Ee
\vskip .1cm
We say the pair of weights $(\sigma,\omega)$ belongs to $\mathcal{C}_{p,q,\alpha}$, if
\Be\label{eq:Cpq}
[\sigma,\omega]_{\mathcal{C}_{p,q,\alpha}}:=\sup_{I\subset \mathbb{R}}|Q_I|_\alpha^{\frac{\gamma}{2+\alpha}-1}\left(\int_{Q_I}\omega dV_\alpha\right)^{1/q}\left(\int_{Q_I}\sigma^{-p'}dV_\alpha\right)^{1/p'}<\infty.
\Ee
\vskip .cm
We say the pair of weights $(\sigma,\omega)$ belongs to $\mathcal{S}_{p,q,\alpha}$, if
\Be\label{eq:Spqclass}[\sigma,\omega]_{\mathcal{S}_{p,q,\alpha}}:=\sup_{I\subset\mathbb{R}}\frac{|Q_I|_{\omega,\alpha}^{p/q}|Q_I|_{\sigma,\alpha}^{p}}{|Q_I|_\alpha^{p(1-\frac{\gamma}{2+\alpha})+1}}\left(\exp\left(\frac{1}{|Q_I|_\alpha}\int_{Q_I}\log \sigma^{-1}dV_\alpha\right)\right)<\infty.\Ee
For corresponding classes in the real case, we refer to \cite{HyPerez, Moen, Sehba1}. 
We have the following norm inequalities.
\btheo\label{thm:Apqestim}
Let $1<p\le q<\infty$, $\alpha>-1$, and $0\le \gamma<2+\alpha$. Let $\sigma,\omega$ be weights. Put $u=\sigma^{-p'}$. Then
\Be\label{eq:ineqApq}\|\mathcal{M}_{\alpha,\gamma} \sigma f\|_{q,\omega,\alpha}\lesssim ([\sigma,\omega]_{\mathcal{A}_{p,q,\alpha}})^{1/p}[\sigma]_{\mathcal{B}_{\infty,\alpha}}^{1/p}\|f\|_{p,\sigma,\alpha},
\Ee
\Be\label{eq:ineqCpq}
\|\mathcal{M}_{\alpha,\gamma} f\|_{q,\omega,\alpha}\lesssim [\sigma,\omega]_{\mathcal{C}_{p,q,\alpha}}[u]_{\mathcal{B}_{\infty,\alpha}}^{\frac{1}{p}}\|\sigma f\|_{p,\alpha},
\Ee
and
\Be\label{eq:ineqSpq}
\|\mathcal{M}_{\alpha,\gamma} \sigma f\|_{q,\omega,\alpha}\lesssim ([\sigma,\omega]_{\mathcal{S}_{p,q,\alpha}})^{1/p}\|f\|_{p,\sigma,\alpha}.
\Ee
\etheo
It is possible to improve (\ref{eq:ineqCpq}) as follows. 
\btheo\label{thm:Cpqestim}
Let $\alpha>-1$,  $0\le \gamma<2+\alpha$, and $1<p\le q<\infty$.  Let $\sigma,\omega$ be weights and put $u=\sigma^{-p'}$. Then

\Be\label{eq:ineqCpqsharp1}
\|\mathcal{M}_{\alpha,\gamma} f\|_{q,\omega,\alpha}\lesssim [\sigma,\omega]_{\mathcal{C}_{p,q,\alpha}}[u]_{\mathcal{B}_{\infty,\alpha}}^{\frac{1}{q}}\|\sigma f\|_{p,\alpha}.
\Ee
\etheo
In particular, when $\sigma=\omega$, writing $$[\omega]_{\mathcal{B}_{p,q,\alpha}}:=\sup_{I\subset \mathbb{R}}\left(\frac{1}{|Q_I|_\alpha}\int_{Q_I}\omega^q dV_\alpha\right)\left(\frac{1}{|Q_I|_\alpha}\int_{Q_I}\omega^{-p'}dV_\alpha\right)^{q/p'},$$
we obtain the following which is sharp.
\bprop\label{prop:Cpqestim}
Let $\alpha>-1$,  $0\le \gamma<2+\alpha$, and $1<p<\frac{2+\alpha}{\gamma}$. Suppose $q$ is defined by the relation $\frac{1}{q}=\frac{1}{p}-\frac{\gamma}{2+\alpha}$.  Let $\omega$ be weights and put $u=\omega^{-p'}$. Then

\Be\label{eq:ineqCpqsharp}
\|\omega \mathcal{M}_{\alpha,\gamma} f\|_{q,\alpha}\lesssim [\omega]_{\mathcal{B}_{p,q,\alpha}}^{\frac{1}{q}}[u]_{\mathcal{B}_{\infty,\alpha}}^{\frac{1}{q}}\|\omega f\|_{p,\alpha},
\Ee
and this is sharp.
\eprop
We also have the following estimate.
\btheo\label{thm:sharpmaxfrac}
Let $\alpha>-1$,  $0\le \gamma<2+\alpha$, and $1<p<\frac{2+\alpha}{\gamma}$. Suppose $q$ is defined by the relation $\frac{1}{q}=\frac{1}{p}-\frac{\gamma}{2+\alpha}$. If $\omega\in \mathcal{B}_{p,q,\alpha}$, then 
\Be\label{eq:ineqBpq}
\|\omega \mathcal{M}_{\alpha,\gamma} f\|_{q,\alpha}\le ([\omega]_{\mathcal{B}_{p,q,\alpha}})^{\frac{p'}{q}(1-\frac{\gamma}{2+\alpha})}\|\omega f\|_{p,\alpha}.
\Ee
Moreover, the exponent $\frac{p'}{q}(1-\frac{\gamma}{2+\alpha})$ is sharp.

\etheo
Taking $p=q$, and putting $\omega=\sigma^{1/p}$, we obtain the following (see also \cite{PR}).
\bcor\label{cor:sharpmax}
Let $\alpha>-1$, and $1<p<\infty$. Then 
\Be\label{eq:ineqsharpmax}
\| \sigma^{1/p}\mathcal{M}_{\alpha} f\|_{p,\sigma,\alpha}\le [\sigma]_{\mathcal{B}_{p,\alpha}}^{\frac{p}{p'}}\|\sigma^{1/p}f\|_{p,\alpha}.
\Ee
\ecor

\vskip .3cm
To prove the sufficient part in the above theorems, we will observe that the matter can be reduced to the case of the dyadic analogue of the corresponding operators. We then appeal to analogues of techmiques of real harmonic analysis:  discritizing  integrals using appropriate level sets, Carleson embeddings and some others (see for example \cite{CarnotBenoit, Cruz1,cuervafrancia,Sawyer, Sehba1} for some of these techniques). We particularly take advantage of the nice properties of the upper-halve of Carleson squares. 
\vskip .1cm
 Given two positive quantities $A$ and $B$, the notation $A\lesssim B$ (resp. $B\lesssim A$) will mean that there is an universal constant $C>0$ such that $A\le CB$ (resp. $B\le CA$). When $A\lesssim B$ and $B\lesssim A$, we write $A\backsimeq B$.
\section{Useful observations and results}
Given an interval $I\subset \mathbb R$, the upper-half of the Carleson box $Q_I$ associated to $I$ is the subset $T_I$ defined by $$T_I:=\{z=x+iy\in \mathbb C: x\in I,\,\,\,\textrm{and}\,\,\,\frac{|I|}{2}<y<|I|\}.$$
Note that $|Q_I|_\alpha\backsimeq |T_I|_\alpha$. 
We consider the following system of dyadic grids,
$$\mathcal D^\beta:=\{2^j\left([0,1)+m+(-1)^j\beta\right):m\in \mathbb Z,\,\,\,j\in \mathbb Z \},\,\,\,\textrm{for}\,\,\,\beta\in \{0,1/3\}.$$
When $\beta=0$, we observe that $\mathcal D^0$ is the standard dyadic grid of $\mathbb R$, denoted $\mathcal D$. 
We recall with \cite{PR} that given an interval $I\subset \mathbb R$, there is a dyadic interval $J\in \mathcal D^\beta$ for some $\beta\in \{0,1/3\}$ such that $I\subseteq J$ and $|J|\le 6|I|$. It follows  that for any locally integrable function $f$,
\Be\label{eq:Maxfunctdyaineq}
\mathcal{M}_{\alpha,\gamma}f(z)\le C\sum_{\beta\in \{0,1/3\}}\mathcal{M}_{\alpha,\gamma}^{d,\beta} f(z),\,\,\,z\in \mathcal H
\Ee 
where $\mathcal{M}_{\alpha,\gamma}^{d,\beta}$ is defined as $\mathcal{M}_{\alpha,\gamma}$ but with the supremum taken only over dyadic intervals of the dyadic grid $\mathcal D^\beta$. 
\vskip .2cm
\begin{definition} Let $\alpha>-1$ and $\omega$ be a positive weight. For any $\delta\geq 1$, a sequence of positive numbers $\{\lambda_{Q_I}\}_{I\in \mathcal{D}^\beta}$ is called a $(\omega,\alpha,\delta)$-Carleson sequence, if there is a constant $C>0$ such that for any $J\in \mathcal{D}^\beta$,
$$\sum_{I\subseteq J,\,I\in \mathcal{D}^\beta}\lambda_{Q_I}\le C|Q_J|_{\omega,\alpha}^\delta.$$
\end{definition}
The smallest constant $C$ in the above definition is called the Carleson constant of the sequence.
\vskip .3cm
It is easy to check that for any $\omega\in \mathcal{B}_{\infty,\alpha}$, the sequence $\{|Q_I|_{\omega,\alpha}\}_{I\in \mathcal{D}^\beta}$ is a $(\omega,\alpha,1)$-Carleson sequence. Indeed, we have
\Beas
\sum_{I\subseteq J,\,I\in \mathcal{D}^\beta}|Q_I|_{\omega,\alpha} &=& \sum_{I\subseteq J,\,I\in \mathcal{D}^\beta}\frac{|Q_I|_{\omega,\alpha}}{|Q_I|_\alpha}|Q_I|_\alpha\\ &\simeq& \sum_{I\subseteq J,\,I\in \mathcal{D}^\beta}\frac{|Q_I|_{\omega,\alpha}}{|Q_I|_\alpha}|T_I|_\alpha\\ &\le& \sum_{I\subseteq J,\,I\in \mathcal{D}^\beta}\int_{T_I}\frac{|Q_I|_{\omega,\alpha}}{|Q_I|_\alpha}dV_\alpha(z)\\ &\le& \sum_{I\subseteq J,\,I\in \mathcal{D}^\beta}\int_{T_I}\mathcal{M}_\alpha(\chi_{Q_J}\omega)dV_\alpha(z)\\ &=& \int_{Q_J}\mathcal{M}_\alpha(\chi_{Q_J}\omega)dV_\alpha(z)\\ &\le& [\omega]_{\mathcal{B}_{\infty,\alpha}}|Q_J|_{\omega,\alpha}.
\Eeas
Let $f\in L^1(\mathcal{H}, dV_\alpha)$ and $\sigma$ a weight. Define the weighted fractional maximal function $\mathcal{M}_{\sigma,\alpha,\gamma}$ as follows:
\Be\label{eq:defweightedfracmax}\mathcal{M}_{\sigma,\alpha,\gamma}f:=\sup_{I\subset \mathbb{R}}\frac{\chi_{Q_I}}{|Q_I|_{\sigma,\alpha}^{1-\frac{\gamma}{2+\alpha}}}\int_{Q_I}|f(z)|\sigma(z) dV_\alpha(z).\Ee
The proof of the following Carleson embedding follows as in \cite{Sehba1}.
\btheo\label{thm:Carlembed0}
Let $\alpha>-1$, and $0\le \gamma<2+\alpha$. Let $\omega$ be a weight on $\mathcal{H}$ and $s\ge 1$. Assume  $\{\lambda_{Q_I}\}_{I\in \mathcal D^\beta}$ is a sequence of positive numbers. If
there exists some constant $A>0$ such that for any interval $J\in \mathcal D^\beta$,
    \begin{equation*}\sum_{I\subseteq J, I\in \mathcal D^\beta}\lambda_{Q_I}\le A|Q_J|_{\omega,\alpha}^s,\end{equation*}
then for any $p\in (1,\infty)$, 
\begin{equation*}\sum_{I, I\in \mathcal D^\beta}\lambda_{Q_I}\left(\frac{1}{|Q_I|_{\omega,\alpha}^{1-\frac{\gamma}{2+\alpha}}}\int_{Q_I}|f(z)|\omega(z)dV_\alpha(z)\right)^{ps}\le A\gamma \left(\int_{\mathcal{H}}(\mathcal{M}_{\omega,\alpha,\gamma}^{d,\beta} f(z))^p\omega(z)dV_\alpha(z)\right)^{s}.\end{equation*}

\etheo

We will also need the following lemma which proof is essentially the same as in the case $\gamma=0$ in \cite{CarnotBenoit}.
\blem\label{lem:equivdefBpq}
Let  $1\le p, q <\infty$ and suppose that $\omega$ is a weight, and $\mu$ a positive measure on $\mathcal H$. Then the following assertions are equivalent.
\begin{itemize}
\item[(i)] There exists a constant $C_1>0$ such that for any interval $I\subset \mathbb R$,
\Be\label{eq:equivMuBpq1}
|Q_I|^{q(\frac{\gamma}{2+\alpha}-\frac{1}{p})}\left(\frac{1}{|Q_I|_\alpha}\int_{Q_I}\omega^{1-p'}(z)dV_\alpha(z)\right)^{q/p'}\mu(Q_I)\le C_1
\Ee
where $\left(\frac{1}{|Q_I|_\alpha}\int_{Q_I}\omega^{1-p'}(z)dV_\alpha(z)\right)^{1/p'}$ is understood as $\left(\inf_{Q_I}\omega\right)^{-1}$ when $p=1$.
\item[(ii)] There exists a constant $C_2>0$ such that for any locally integrable function $f$ and  any interval $I\subset \mathbb R$,
\Be\label{eq:equivMuBpq2}
\left(\frac{1}{|Q_I|_\alpha^{1-\frac{\gamma}{2+\alpha}}}\int_{Q_I}|f(z)|dV_\alpha(z)\right)^q\mu(Q_I)\le C_2\left(\int_{Q_I}|f(z)|^p\omega(z) dV_\alpha(z)\right)^{q/p}.
\Ee
\end{itemize}
\elem

\section{Proof of the results}
\subsection{Proof of Theorem \ref{thm:main1}}
We start with the following level sets embedding. The proof follows exactly as in case $\gamma=0$ (see \cite[Lemma 3.4.]{CarnotBenoit})
\blem\label{lem:levelsetsembed}
Let $f$ be a locally integrable function. Then for any $\lambda>0$,
\Be\label{eq:levelsetsembed}
\{z\in \mathcal H: \mathcal{M}_{\alpha,\gamma}f(z)>\lambda\}\subset \{z\in \mathcal{H}: \mathcal{M}_{\alpha,\gamma}^{d,\beta}f(z)>\frac{\lambda}{C_{\alpha,\gamma}}\}
\Ee
where $C_{\alpha,\gamma}=2^{2+\alpha-\gamma}(1+2^{4+2\alpha-2\gamma})$.
\elem
\begin{proof}[Proof of Theorem \ref{thm:main1}]
Let us start by proving (b). Recall that when $\gamma=0$, $\mathcal{M}_{\alpha,\gamma}=\mathcal{M}_{\alpha}$, and $\|\mathcal{M}_{\alpha}\|_{\infty}\le \|f\|_\infty$. Next assume that $\gamma\neq 0$. For $z\in \mathcal{H}$, let $Q_I$ be a Carleson square containing $z$. Then
\Beas
\frac{1}{|I|^{2+\alpha-\gamma}}\int_{Q_I}|f(w)|dV_\alpha(w) &\le& \frac{1}{|Q_I|_\alpha^{1-\frac{\gamma}{2+\alpha}}}\|f\|_{\frac{2+\alpha}{\gamma},\alpha}\left(\int_{Q_I}dV_\alpha(w)\right)^{1-\frac{\gamma}{2+\alpha}}\\ &=& \frac{1}{|Q_I|_\alpha^{1-\frac{\gamma}{2+\alpha}}}\|f\|_{\frac{2+\alpha}{\gamma},\alpha}|Q_I|_\alpha^{1-\frac{\gamma}{2+\alpha}}\\ &=& \|f\|_{\frac{2+\alpha}{\gamma},\alpha}.
\Eeas
To prove the inequality (\ref{eq:main11}), it is enough by Lemma \ref{lem:levelsetsembed} to prove the same inequality with $\mathcal{M}_{\alpha,\gamma}$ replaced by $\mathcal{M}_{\alpha,\gamma}^{d,\beta}$. 
It is then enough to prove the following.
\bprop\label{prop:weightedfracmax}
Let $\alpha>-1$, and $0\le \gamma<2+\alpha$. Assume $\sigma$ is a weight. Then for any $\lambda>0$, there exists a positive constant $C$ such that
\Be\label{eq:weightedfracmax}|\{z\in \mathcal{H}:\mathcal{M}_{\sigma,\alpha,\gamma}^{d,\beta}f(z)>\lambda\}|_{\sigma,\alpha}\le C\left(\frac{1}{\lambda}\int_{\mathcal{H}}|f(z)|\sigma(z)dV_\alpha(z)\right)^{\frac{2+\alpha}{2+\alpha-\gamma}}.
\Ee
\eprop
\begin{proof}
Put $$E_\lambda^{d,\beta}:=\{z\in \mathcal{H}: \mathcal{M}_{\sigma,\alpha,\gamma}^{d,\beta}f(z)>\lambda\}.$$ Then following usual arguments, $E_\lambda^{d,\beta}=\bigcup_I Q_I$ where $(Q_I)_{I\in \mathcal{D}^\beta}$ is a family of maximal Carleson squares. In particular, if $z\in E_\lambda^{d,\beta}$, then there is an interval $I\subset \mathcal{D}^\beta$ such that $z\in Q_I$ and
$$\frac{1}{|Q_I|_{\sigma,\alpha}^{1-\frac{\gamma}{2+\alpha}}}\int_{Q_I}|f(w)|\sigma(w)dV_\alpha(w)>\lambda,$$ 
that is $$|Q_I|_{\sigma,\alpha}\le \left(\frac{1}{\lambda}\int_{Q_I}|f(w)|\sigma(w)dV_\alpha(w)\right)^{\frac{2+\alpha}{2+\alpha-\gamma}}.$$
Hence
\Beas
M &:=& |\{z\in \mathcal H: \mathcal{M}_{\sigma,\alpha,\gamma}^df(z)>\lambda\}|_{\sigma,\alpha}\\ &\le& |E_\lambda^{d,\beta}|_{\sigma,\alpha}\le \sum_I|Q_I|_{\sigma,\alpha}\\ &\le& \sum_I \left(\frac{1}{\lambda}\int_{Q_I}|f(w)|\sigma(w)dV_\alpha(w)\right)^{\frac{2+\alpha}{2+\alpha-\gamma}}\\ &\le& \left(\frac{1}{\lambda}\sum_I\int_{Q_I}|f(w)|\sigma(w)dV_\alpha(w)\right)^{\frac{2+\alpha}{2+\alpha-\gamma}}\\ &\le& \left(\frac{1}{\lambda}\int_{\{z\in \mathcal{H}: \mathcal{M}_{\sigma,\alpha,\gamma}^{d,\beta}f(z)>\lambda\}}|f(w)|\sigma(w)dV_\alpha(w)\right)^{\frac{2+\alpha}{2+\alpha-\gamma}}\\ &\le& \left(\frac{1}{\lambda}\int_{\mathcal{H}}|f(w)|\sigma(w)dV_\alpha(w)\right)^{\frac{2+\alpha}{2+\alpha-\gamma}}
\Eeas
\end{proof}
The proof is complete.
\end{proof}
We observe that from the above discussion and Marcinkiewicz interpolation theorem one has the following useful result.
\bcor\label{cor:weightedfracmax}
Let $\alpha>-1$, $0\le \gamma<2+\alpha$, and let $\sigma$ be a weight. Assume that $q$ is such that $\frac{1}{p}-\frac{1}{q}=\frac{\gamma}{2+\alpha}$. Then there exists a constant $C=C(p,\alpha,\gamma)$ such that 
\Be\label{eq:offdiagweightedfracmax}
\left(\int_{\mathcal{H}}(\mathcal{M}_{\sigma,\alpha,\gamma}^{d,\beta}f(z))^q\sigma(z)dV_\alpha(z)\right)^{1/q}\le C\left(\int_{\mathcal{H}}|f(z)|^p\sigma(z)dV_\alpha(z)\right)^{1/p}
\Ee
 
\ecor
We specify the constant for $\mathcal{M}_{\alpha,\gamma}$.
\bprop
Let $\alpha>-1$. If $0\le \gamma<2+\alpha$, $1<p<\frac{2+\alpha}{\gamma}$, and $\frac{1}{q}=\frac{1}{p}-\frac{\gamma}{2+\alpha}$, then
\Be\label{eq:offdiagunweightedfracmax}
\left(\int_{\mathcal{H}}(\mathcal{M}_{\alpha,\gamma}f(z))^qdV_\alpha(z)\right)^{1/q}\le \left((1+\frac{p'}{q})C_{\alpha,\gamma}\right)^{1-\frac{\gamma}{2+\alpha}}\left(\int_{\mathcal{H}}|f(z)|^pdV_\alpha(z)\right)^{1/p}
\Ee
where $C_{\alpha,\gamma}$ is the constant in (\ref{eq:levelsetsembed}).
\eprop
\begin{proof}
From Lemma \ref{lem:levelsetsembed} and the proof of Theorem \ref{thm:main1}, we have
$$|\{z\in \mathcal H: \mathcal{M}_{\alpha,\gamma}f(z)>\lambda\}|_\alpha\le \left(\frac{1}{\lambda}\int_{\{z\in \mathbb D: \mathcal{M}_{\alpha,\gamma}^{d,\beta}f(z)>\frac{\lambda}{C_{\alpha,\gamma}}\}}|f(w)|dV_\alpha(w)\right)^{\frac{2+\alpha}{2+\alpha-\gamma}}.$$
Thus 
\Beas
L &:=& \int_{\mathcal{H}}(\mathcal{M}_{\alpha,\gamma}f(z))^qdV_\alpha(z)\\ &=& \int_0^\infty q\lambda^{q-1}|\{z\in \mathcal H: \mathcal{M}_{\alpha,\gamma}f(z)>\lambda\}|_\alpha d\lambda\\ &\le& \int_0^\infty q\lambda^{q-1}\left(\frac{1}{\lambda}\int_{\{z\in \mathbb D: \mathcal{M}_{\alpha,\gamma}^{d,\beta}f(z)>\frac{\lambda}{C_{\alpha,\gamma}}\}}|f(w)|dV_\alpha(w)\right)^{\frac{2+\alpha}{2+\alpha-\gamma}}d\lambda\\ &\le& \left(\int_{\mathcal{H}}|f(z)|\left(\int_0^{C_{\alpha,\gamma}\mathcal{M}_{\alpha,\gamma}f(z)}\lambda^{q-\frac{2+\alpha}{2+\alpha-\gamma}-1}d\lambda\right)^{\frac{2+\alpha-\gamma}{2+\alpha}}dV_\alpha(z)\right)^{\frac{2+\alpha}{2+\alpha-\gamma}}\\ &\le& \frac{qC_{\alpha,\gamma}}{q-\frac{2+\alpha}{2+\alpha-\gamma}}\left(\int_{\mathcal{H}}|f(z)|(\mathcal{M}_{\alpha,\gamma}f(z))^{\frac{q}{p'}}dV_\alpha(z)\right)^{\frac{2+\alpha}{2+\alpha-\gamma}}\\ &\le& \frac{qC_{\alpha,\gamma}}{q-\frac{2+\alpha}{2+\alpha-\gamma}}\left(\int_{\mathcal{H}}|f(z)|^pdV_\alpha(z)\right)^{\frac{2+\alpha}{p(2+\alpha-\gamma)}}\left(\int_{\mathcal{H}}(\mathcal{M}_{\alpha,\gamma}f(z))^qdV_\alpha(z)\right)^{\frac{2+\alpha}{2+\alpha-\gamma}\frac{1}{p'}}.
\Eeas
That is 
$$\left(\int_{\mathcal{H}}(\mathcal{M}_{\alpha,\gamma}f(z))^qdV_\alpha(z)\right)^{1-\frac{2+\alpha}{2+\alpha-\gamma}\frac{1}{p'}}\le \frac{qC_{\alpha,\gamma}}{q-\frac{2+\alpha}{2+\alpha-\gamma}}\left(\int_{\mathcal{H}}|f(z)|^pdV_\alpha(z)\right)^{\frac{2+\alpha}{p(2+\alpha-\gamma)}},$$
which is equivalent to
$$\left(\int_{\mathcal{H}}(\mathcal{M}_{\alpha,\gamma}f(z))^qdV_\alpha(z)\right)^{1/q}\le \left(\frac{p'+q}{q}C_{\alpha,\gamma}\right)^{\frac{1}{p'}+\frac{1}{q}}\left(\int_{\mathcal{H}}|f(z)|^pdV_\alpha(z)\right)^{1/p}.$$
The proof is complete.
\end{proof}
\subsection{Proof of Theorem \ref{thm:main2}}
Let us prove Theorem \ref{thm:main2}.
\begin{proof}[Proof of Theorem~ {\rm\ref{thm:main2}}]
Let us note that by Lemma \ref{lem:equivdefBpq}, $(\textrm{b})\Leftrightarrow (\textrm{c})$. Let us prove that $(\textrm{a})\Leftrightarrow (\textrm{c})$.

Let $f$ be a locally integrable function and $I$ an interval. Fix $\lambda$ such that $$0<\lambda<\frac{1}{|Q_I|_\alpha^{1-\frac{\gamma}{2+\alpha}}}\int_{Q_I}|f(z)|dV_\alpha(z).$$ Then \[Q_I\subset \{z\in \mathcal H:\mathcal{M}_{\alpha,\gamma}(\chi_{Q_I}f)>\lambda)\}.\]
It follows from the latter and (\ref{eq:main21}) that \[\mu(Q_I)\le \frac{C}{\lambda^q}\left(\int_{Q_I}|f(z)|^p\omega(z)dV(z)\right)^{q/p}.\]
As this happens for all $\lambda>0$, it follows in particular that
\[\mu(Q_I)\left(\frac{1}{|Q_I|_\alpha^{1-\frac{\gamma}{2+\alpha}}}\int_{Q_I}|f(z)|dV_\alpha(z)\right)^q\le C\left(\int_{Q_I}|f(z)|^p\omega(z)dV(z)\right)^{q/p}.\]
\vskip .1cm
Next suppose that (\ref{eq:main23}) holds. We observe with Lemma \ref{lem:levelsetsembed} that to obtain (\ref{eq:main21}), we only have to prove the following
\Be\label{eq:main34}
\mu\left(\{z\in \mathcal D: \mathcal{M}_{\alpha,\gamma}^{d,\beta}f(z)>\frac{\lambda}{C_{\alpha,\gamma}}\}\right)\le \frac{C}{\lambda^q}\|f\|_{p,\omega,\alpha}^q.
\Ee
We recall that \[\{z\in \mathcal H: \mathcal{M}_{\alpha,\gamma}^{d,\beta}f(z)>\frac{\lambda}{C_{\alpha,\gamma}}\}=\cup_{j\in \mathbb N_0}Q_{I_j}\]
where the $I_j$s are maximal dyadic intervals (in $\mathcal{D}^\beta$) with respect to the inclusion and such that
\[\frac{1}{|Q_{I_j}|_\alpha^{1-\frac{\gamma}{2+\alpha}}}\int_{Q_{I_j}}|f|dV_\alpha>\frac{\lambda}{C_{\alpha,\gamma}}.
\]
Our hypothesis provides in particular that
\[
\mu(Q_{I_j})\le C\left(\frac{|Q_I|_\alpha^{1-\frac{\gamma}{2+\alpha}}}{\int_{Q_{I_j}}|f|dV_\alpha}\right)^q\left(\int_{Q_{I_j}}|f|^p\omega dV_\alpha\right)^{q/p}.
\]

Thus
\begin{eqnarray*}
\mu\left(\{z\in \mathcal {H}: \mathcal{M}_{\alpha,\gamma}^{d,\beta}f(z)>\frac{\lambda}{C_{\alpha,\gamma}}\}\right) &=& \sum_{j}\mu(Q_{I_j})\\ &\le& C\sum_{j}\left(\frac{|Q_I|_\alpha^{1-\frac{\gamma}{2+\alpha}}}{\int_{Q_{I_j}}|f|dV_\alpha}\right)^q\left(\int_{Q_{I_j}}|f|^p\omega dV_\alpha\right)^{q/p}\\ &\le& \left(\frac{C_{\alpha,\gamma}}{\lambda}\right)^q\sum_{j}\left(\int_{Q_{I_j}}|f|^p\omega dV_\alpha\right)^{q/p}\\ &\le& \left(\frac{C_{\alpha,\gamma}}{\lambda}\right)^q\left(\sum_{j}\int_{Q_{I_j}}|f|^p\omega dV_\alpha\right)^{q/p}\\ &\le& \left(\frac{C_{\alpha,\gamma}}{\lambda}\right)^q\left(\int_{\mathcal H}|f(z)|^p\omega(z) dV_\alpha(z)\right)^{q/p}\\ &=& \left(\frac{C_{\alpha,\gamma}}{\lambda}\right)^q\|f\|_{p,\omega,\alpha}^q.
\end{eqnarray*}
The proof is complete.
\end{proof}
Taking $d\mu(z)=\sigma(z)dV(z)$, we obtain the following corollary.
\begin{corollary}\label{cor:main4cor}
Let $\alpha>-1$, $0\le \gamma<2+\alpha$, and $1\le p\le q<\infty$. Let $\omega,\sigma$ be  two weights on $\mathcal H$. Then the following assertions are equivalent.
\begin{itemize}
\item[(a)] There is a constant $C_1>0$ such that for any $f\in L_\omega^p(\mathcal H)$, and any $\lambda>0$,
\begin{equation}\label{eq:main41cor}
|\{z\in \mathcal {H}: \mathcal{M}_{\alpha,\gamma}f(z)>\lambda\}|_{\sigma,\alpha}\le \frac{C_1}{\lambda^q}\left(\int_{\mathcal H}|f(z)|^p\omega(z)dV_\alpha(z)\right)^{q/p}
\end{equation}
\item[(b)] There is a constant $C_2>0$ such that for any interval $I\subset \mathbb {R}$,
\begin{equation}\label{eq:main42cor}
|Q_I|_\alpha^{\frac{\gamma}{2+\alpha}+1/q-1/p}\left(\frac{1}{|Q_I|_\alpha}\int_{Q_I}\omega^{1-p'}(z)dV_\alpha(z)\right)^{1/p'}\left(\frac{1}{|Q_I|_\alpha}\int_{Q_I}\sigma(z)dV_\alpha(z)\right)^{1/q}\le C_2
\end{equation}
where $\left(\frac{1}{|Q_I|_\alpha}\int_{Q_I}\omega^{1-p'}(z)dV_\alpha(z)\right)^{1/p'}$ is understood as $\left(\inf_{Q_I}\omega\right)^{-1}$ when $p=1$.
\end{itemize}
\end{corollary}
\subsection{Proof of Theorem \ref{thm:main3}}
The proof of Sawyer-type characterization for the maximal functions is a routine. We adopt the classical proof here (see for example \cite{Sehba1}).
\begin{proof}[Proof of Theorem \ref{thm:main3}]
That (a)$\Rightarrow$ (b) follows by taking $f=\chi_{Q_I}$ in (\ref{eq:main31}). To prove that (b)$\Rightarrow$ (a), it is enough by the observations made in the last section, to prove that under (\ref{eq:main32}), for any $\beta=0,\frac{1}{3}$, there is a positive constant $C$ such that for any function $f$,
\Be\label{eq:main33dya}
\int_{\mathcal{H}}(\mathcal{M}_{\alpha,\gamma}^{d,\beta}(\sigma f)(z))^qd\mu(z)\le C\left(\int_{\mathcal{H}}|f(z)|^p\sigma(z)dV_\alpha(z)\right)^{q/p}
\Ee
Let us prove (\ref{eq:main33dya}): let $a\geq 2^{2+\alpha-\gamma}$. To each integer $k$, we associate the set $$\Omega_k:=\{z\in \mathcal{H}: a^k<\mathcal{M}_{\alpha,\gamma}^{d,\beta}(\sigma f)(z)\le a^{k+1}\}.$$
Then we have that $\Omega_k\subset \bigcup_{j=1}^\infty Q_{I_j}^k$ where $I_j\in \mathcal{D}^\beta$ and $\{Q_{I_j}^k\}_{j\geq 1}$ is a family of dyadic Carleson square which is maximal with respect to the inclusion and such that $$2^{2+\alpha-\gamma}a^k>\frac{1}{|Q_{I_j}^k|_\alpha^{1-\frac{\gamma}{2+\alpha}}}\int_{Q_{I_j}^k}|f(z)|\sigma(z)dV_\alpha(z)>a^k.$$
Define $E(Q_{I_j}^k):=Q_{I_j}^k\cap \Omega_k$. Then $\Omega_k=\bigcup_{j=1}^\infty E_{k,j}$ and the $E_{k,j}$ are disjoint for all $k$ and $j$, that is $E_{k,j}\cap E_{l,m}=\emptyset$ for $(k,j)\neq (l,m)$. It follows that
\Beas
L &:=& \int_{\mathbb R^n}\left(\mathcal M_{\alpha,\gamma}^{d,\beta}({\sigma}{f})(z)\right)^qd\mu(z)\\ &=& \sum_k\int_{\Omega_k}\left(\mathcal M_{\alpha,\gamma}^{d,\beta}({\sigma}{f})(z)\right)^qd\mu(z)\\ &\le&
a^q\sum_ka^{kq}\mu(\Omega_k)\\ &\le& a^q\sum_{k,j}a^{kq}\mu(E(Q_{I_j}^k))\\ &\le& a^q\sum_{k,j}\left(\frac{1}{|Q_{I_j}^k|_\alpha^{1-\frac{\gamma}{2+\alpha}}}\int_{Q_{I_j}^k}|f(z)|\sigma(z)dV_\alpha(z)\right)^q\mu(E(Q_{I_j}^k))\\ &=& a^q\sum_{k,j}\left(\frac{1}{|Q_{I_j}^k|_{\sigma,\alpha}}\int_{Q_{I_j}^k}|f(z)|\sigma(z)dV_\alpha(z)\right)^q\frac{\mu(E(Q_{I_j}^k))|Q_{I_j}^k|_{\sigma,\alpha}^q}{|Q_{I_j}^k|_\alpha^{q(1-\frac{\gamma}{2+\alpha})}}\\ &\le& C\left(\int_{\mathcal{H}}|f(z)|^p\sigma(z)dV_\alpha(z)\right)^{q/p}
\Eeas
provided the sequence
$$
\lambda_{Q_I}:=\left\{ \begin{matrix} \frac{\mu(E(Q_I))|Q_I|_{\sigma,\alpha}^q}{|Q_I|_\alpha^{q(1-\frac{\gamma}{2+\alpha})}} &\text{if }& Q_I=Q_{I_j}^k\,\,\,\textrm{for some}\,\,\,(k,j),\\
      0 & \text{ otherwise} .
                                  \end{matrix} \right.
$$
is a $(\sigma,\alpha,\frac{q}{p})$-Carleson sequence. Using (\ref{eq:main32}),we obtain for any $J\in \mathcal{D}^\beta$
\Beas
\sum_{I\subseteq J,\, I\in \mathcal{D}^\beta }\lambda_{Q_I} &=& \sum_{I\subseteq J,\, I\in \mathcal{D}^\beta }\frac{\mu(E(Q_I))|Q_I|_{\sigma,\alpha}^q}{|Q_I|_\alpha^{q(1-\frac{\gamma}{2+\alpha})}}\\ &=& \sum_{I\subseteq J,\, I\in \mathcal{D}^\beta }\left(\frac{1}{|Q_I|_\alpha^{1-\frac{\gamma}{2+\alpha}}}|Q_I|_{\sigma,\alpha}\right)^q\mu(E(Q_I))\\ &\leq& \sum_{I\subseteq J,\, I\in \mathcal{D}^\beta }\left(\frac{1}{|Q_I|_\alpha^{1-\frac{\gamma}{2+\alpha}}}\int_{Q_I}\chi_{Q_J}(z)\sigma(z) dV_\alpha(z)\right)^q\mu(E(Q_I))\\ &\leq& \sum_{I\subseteq J,\, I\in \mathcal{D}^\beta }\int_{E(Q_I)}(\mathcal{M}_{\alpha,\gamma}^{d,\beta}(\sigma\chi_{Q_J})(z))^qd\mu(z)\\ &=& \int_{Q_J}(\mathcal{M}_{\alpha,\gamma}^{d,\beta}(\sigma\chi_{Q_J})(z))^qd\mu(z)\\ &\le&  C|Q_J|_{\sigma,\alpha}^{q/p}.
\Eeas
That is $\{\lambda_{Q_I}\}_{i\in \mathcal{D}^\beta}$ is a $(\sigma,\alpha,\frac{q}{p})$-Carleson sequence. The proof is complete.
\end{proof} 
\subsection{Proof of Theorem \ref{thm:main4}}
First suppose that (\ref{eq:main41}) holds and observe that for any interval $I\subset \mathbb R$, $\frac{|Q_I|_{\sigma,\alpha}}{|Q_I|_\alpha^{1-\frac{\gamma}{2+\alpha}}}\le \mathcal{M}_{\alpha,\gamma} (\sigma\chi_{Q_I})(z)$ for any $z\in Q_I$. It follows that
\Beas\left(\frac{|Q_I|_{\sigma,\alpha}^q}{|Q_I|_\alpha^{q(1-\frac{\gamma}{2+\alpha})}}\mu(Q_I)\right)^{1/q} &\le& \left(\int_{\mathcal{H}}\left(\mathcal{M}_{\alpha,\gamma} (\sigma\chi_{Q_I})(z)\right)^qd\mu(z)\right)^{1/q}\\ &\le& C_1\|\chi_{Q_I}\|_{p,\sigma,\alpha}=|Q_I|_{\sigma,\alpha}^{1/p}\Eeas
which provides that for any interval $I\subset \mathbb R$,
$$|Q_I|_\alpha^{-q(1-\frac{\gamma}{2+\alpha})}\mu(Q_I)|Q_I|_{\sigma,\alpha}^{\frac{q}{p'}}\le C_1.$$
That is (\ref{eq:main42}) holds.

To prove that $(\textrm{ii})\Rightarrow (\textrm{i})$, it is enough by the observations made at the beginning of the previous section to prove the following.
\blem\label{lem:main1}
Let $\alpha>-1$, $0\le \gamma<2+\alpha$, and $1< p\le q<\infty$. Assume that $\mu$ is positive Borel measure on $\mathcal{H}$ and $\sigma$ is a weight in the class $\mathcal{B}_{\infty,\alpha}$ such that (\ref{eq:main42}) holds. Then there is a positive constant $C$ such that for any function $f$, and any $\beta\in \{0,\frac{1}{3}\}$,
\Be\label{eq:main13}
\left(\int_{\mathcal H}(\mathcal{M}_{\alpha,\gamma}^{d,\beta} (\sigma f)(z))^qd\mu(z)\right)^{1/q}\le C\left(\int_{\mathcal{H}}|f(z)|^p\sigma(z)dV_\alpha(z)\right)^{1/p}.
\Ee
\elem
\begin{proof}
Let $a\ge 2^{2+\alpha-\gamma}$. To each integer $k$, we associate the set
$$\Omega_{k}:=\{z\in \mathcal H: a^k<\mathcal{M}_{\alpha,\gamma}^{d,\beta}(\sigma f)(z)\leq a^{k+1}\}.$$
We observe that
$\Omega_{k}\subset \cup_{j=1}^{\infty}Q_{I_{j}}^k,$ where
$Q_{I_{j}}^k$ ($I_{j}\in \mathcal{D}^{\beta}$) is a dyadic cube maximal (with respect to the inclusion) such that
$$\frac{1}{|Q_{I_{j}}^k|_\alpha^{1-\frac{\gamma}{2+\alpha}}}\int_{Q_{I_{j}}^k}|(\sigma f)(z)|dV_\alpha(z)>a^k.$$
Following the same reasoning as in the proof of the inequality (\ref{eq:main33dya}), we obtain
 \Beas
L &:=& \int_{\mathcal H}(\mathcal{M}_{\alpha,\gamma}^{d,\beta}(\sigma f)(z))^{q}d\mu(z)\\ &=& \sum_{k}\int_{\Omega_k}(\mathcal{M}_{\alpha,\gamma}^{d,\beta}(\sigma f)(z))^{q}d\mu(z)\\
 &\le&
 a^{q}\sum_{k,j}\left(\frac{1}{|Q_{I_{j}}^k|_\alpha^{1-\frac{\gamma}{2+\alpha}}}\int_{Q_{I_{j}}^k}|(\sigma f)(z)|dV_\alpha(z)\right)^{q}\mu(Q_{I_{j}}^k)\\
 &\lesssim&
 a^{q}\sum_{k,j}\left(\frac{1}{|Q_{I_{j}}^k|_{\sigma,\alpha}}\int_{Q_{I_{j}}^k}|(\sigma f)(z)|dV_\alpha(z)\right)^{q}|Q_{I_j}^k|_\alpha^{-q(1-\frac{\gamma}{2+\alpha})}\mu(Q_{I_j}^k)|Q_{I_j}^k|_{\sigma,\alpha}^{q}\\
&\lesssim& [\sigma,\mu]_{p,q,\alpha,\gamma}a^{q}\sum_{k,j}\left(\frac{1}{|Q_{I_{j}}^k|_{\sigma,\alpha}}\int_{Q_{I_{j}}^k}|(\sigma f)(z)|dV_\alpha(z)\right)^{q}|Q_{I_j}^k|_{\sigma,\alpha}^{q/p}
\\
 &\le& C[\sigma,\mu]_{p,q,\alpha,\gamma}\left(\int_{\mathcal{H}}|f(z)|^p\sigma(z)dV_\alpha(z)\right)^{q/p}
\Eeas
provided the sequence 
$$
\lambda_{Q_I}:=\left\{ \begin{matrix} |Q_{I_j}^k|_{\sigma,\alpha}^{q/p} &\text{if }& Q_I=Q_{I_j}^k\,\,\,\textrm{for some}\,\,\,(k,j),\\
      0 & \text{ otherwise} .
                                  \end{matrix} \right.
$$
is a $(\sigma,\alpha,\frac{q}{p})$-Carleson sequence. We have seen in the previous section that for $\sigma\in \mathcal{B}_{\infty,\alpha}$, $\{|Q_I|_{\sigma,\alpha}\}_{I\in \mathcal{D}^\beta}$ was a $(\sigma,\alpha,1)$-Carleson sequence with Carleson constant $[\omega]_{\mathcal{B}_{\infty,\alpha}}$. Thus $\{|Q_I|_{\sigma,\alpha}^{q/p}\}_{I\in \mathcal{D}^\beta}$ is a $(\sigma,\alpha,\frac{q}{p})$-Carleson sequence with constant $[\omega]_{\mathcal{B}_{\infty,\alpha}}^{q/p}$. The proof is complete.

\end{proof}
Taking $d\mu(z)=\omega(z)dV_\alpha(z)$ where $\omega$ is a weight, we obtain the following corollary.
\bcor\label{cor:main44}
Let $1<p\le q<\infty$, and $\omega, \sigma$ be two weights on $\mathcal H$. Assume that $\sigma\in \mathcal{B}_{\infty,\alpha}$. Then the following assertions are equivalent.
 \begin{itemize}
\item[(i)] There exists a constant $C_1>0$ such that for any $f\in L_\omega^p(\mathcal H)$,
\Be
\left(\int_{\mathcal H}\left(\mathcal{M}_{\alpha,\gamma} (\sigma f)(z)\right)^q\omega(z)dV_\alpha(z)\right)^{1/q}\le C_1\|f\|_{p,\sigma,\alpha}.
\Ee

\item[(ii)] There is a constant $C_2>0$ such that for any interval $I\subset \mathbb {R}$,
\begin{equation}\label{eq:cormain44}
\frac{|Q_I|_{\omega,\alpha}|Q_I|_{\sigma,\alpha}^{q/p'}}{|Q_I|_\alpha^{q(1-\frac{\gamma}{2+\alpha})}}\le C_2
\end{equation}
\end{itemize}
Moreover, $$\|\mathcal{M}_{\alpha,\gamma} \sigma f\|_{q,\omega,\alpha}\le ([\sigma,\omega]_{\mathcal{A}_{p,q,\alpha}})^{1/p}[\sigma]_{\mathcal{B}_{\infty,\alpha}}^{1/p}\|f\|_{p,\sigma,\alpha}.$$
\ecor

\subsection{Proof of Theorem \ref{thm:main5}}
 Let $\Psi$ be the complementary function of the Young function $\Phi$. Recall the following generalized H\"olders's inequality:
\Be\label{eq:holderhgene}
\frac{1}{|Q_{I}|_\alpha}\int_{Q_{I}}|(fg)(z)|dV_\alpha(z)\le \|f\|_{Q_I,\Phi,\alpha}\|g\|_{Q_I,\Psi,\alpha}
\Ee
We recall the following result (see \cite{Sehba3}).
\blem\label{lem:phimax}
Let $\alpha>-1$, $1<p<\infty$.. Then there is a positive constant $C$ such that for any function $f$,

\Be\label{eq:phimax}
\left(\int_{\mathcal H}\left(\mathcal{M}_{\Phi,\alpha} f(z)\right)^p d\mu(z)\right)^{1/p}\le C\left(\int_{\mathcal H}|f(z)|^pdV_\alpha(z)\right)^{1/p}
\Ee
\elem
\begin{proof}[Proof of Theorem \ref{thm:main5}]
As above, to establish the inequality (\ref{eq:main52}), we only need to prove that the same inequality holds with $\mathcal{M}_{\alpha,\gamma}$ replaced by the dyadic maximal function $\mathcal{M}_{\alpha,\gamma}^{d,\beta}$, $\beta\in \{0, 1/3\}$. Once again, for each integer $k$, we define the set
$$\Omega_{k}:=\{z\in \mathcal H: a^k<\mathcal{M}_{\alpha,\gamma}^{d,\beta}(f)(z)\leq a^{k+1}\}.$$
We already know that
$\Omega_{k}\subset \cup_{j=1}^{\infty}Q_{I_{j}}^k,$ where
$Q_{I_{j}}^k$ ($I_{j}\in \mathcal{D}^{\beta}$) is a dyadic cube maximal (with respect to the inclusion) such that
$$\frac{1}{|Q_{I_{j}}^k|_\alpha^{1-\frac{\gamma}{2+\alpha}}}\int_{Q_{I_{j}}^k}|f(z)|dV_\alpha(z)>a^k.$$
Following the same reasoning as in the proof of Theorem \ref{thm:main4} , we first obtain 
 \Beas
L &:=& \int_{\mathcal H}(\mathcal{M}_{\alpha,\gamma}^{d,\beta} f(z))^{q}d\mu(z)\\ 
 &\le&
 a^{q}\sum_{k,j}\left(\frac{1}{|Q_{I_{j}}^k|_\alpha^{1-\frac{\gamma}{2+\alpha}}}\int_{Q_{I_{j}}^k}|f(z)|dV_\alpha(z)\right)^{q}\mu(Q_{I_{j}}^k)\\
 &\lesssim&
 a^{q}\sum_{k,j}\left(\frac{1}{|Q_{I_{j}}^k|_{\alpha}}\int_{Q_{I_{j}}^k}|f(z)|dV_\alpha(z)\right)^{q}|Q_{I_j}^k|_\alpha^{\frac{q\gamma}{2+\alpha}}\mu(Q_{I_j}^k).
\Eeas
Now using (\ref{eq:holderhgene}) and (\ref{eq:main51}), we obtain
\Beas 
L &\lesssim& a^{q}\sum_{k,j}\left(\frac{1}{|Q_{I_{j}}^k|_{\alpha}}\int_{Q_{I_{j}}^k}|(\omega f)(z)|\omega^{-1}(z)dV_\alpha(z)\right)^{q}|Q_{I_j}^k|_\alpha^{\frac{q\gamma}{2+\alpha}}\mu(Q_{I_j}^k)
\\ &\le& a^{q}\sum_{k,j}\|f\omega\|_{Q_{I_j}^k,\Phi,\alpha}^q\|\omega^{-1}\|_{Q_{I_j}^k,\Psi,\alpha}^q|Q_{I_j}^k|_\alpha^{\frac{q\gamma}{2+\alpha}}\mu(Q_{I_j}^k)\\ &\le& C\sum_{k,j}\|f\omega\|_{Q_{I_j}^k,\Phi,\alpha}^q|Q_{I_j}^k|_\alpha^{\frac{q}{p}}.
\Eeas
Finally, using Lemma \ref{lem:phimax} and writing $T_{I_j}^k$ for the upper-half of the square $Q_{I_j}^k$, we obtain
\Beas
L &\le& C\left(\sum_{k,j}\|f\omega\|_{Q_{I_j}^k,\Phi,\alpha}^p|Q_{I_j}^k|_\alpha\right)^{\frac{q}{p}}\\  &\le& C\left(\sum_{k,j}\int_{T_{I_j}^k}\|f\omega\|_{Q_{I_j}^k,\Phi,\alpha}^pdV_\alpha(z)\right)^{\frac{q}{p}}\\ &\le& C\left(\int_{\mathcal{H}}(\mathcal{M}_{\Phi,\alpha}(\omega f)(z))^pdV_\alpha(z)\right)^{\frac{q}{p}}\\
 &\le& C\left(\int_{\mathcal{H}}|(\omega f)(z)|^pdV_\alpha(z)\right)^{q/p}.
\Eeas
The proof is complete.
\end{proof}
It is easy to see that for $1<p<\infty$, and $r>1$, $\Phi(t)=t^{(p'r)'}$ is in the class $B_p$. Thus we derive the following.
\bcor\label{cor:main1}
Let $\alpha>-1$, $0\le \gamma<2+\alpha$, and $1<p\le q<\infty$.  Assume that $\omega$ a weight and $\mu$ a positive Borel measure on $\mathcal H$ such that for some $r>1$, there is positive constant $C$ for which for any interval $I\subset \mathbb{R}$
\Be\label{eq:main31off}
|Q_I|_\alpha^{q(\frac{\gamma}{2+\alpha}-\frac{1}{p})}\left(\frac{1}{|Q_I|_\alpha}\int_{Q_I}\omega^{-p'r}dV_\alpha\right)^{q/p'r}\mu(Q_I)\le C.
\Ee

 Then there is a positive constant $K$ such that for any $f\in L^p(\mathcal H, \omega dV_\alpha)$, 
 \Be\label{eq:main32off}
 \left(\int_{\mathcal H}\left(\mathcal{M}_{\alpha,\gamma} f(z)\right)^qd\mu(z)\right)^{1/q}\le K\|f\omega\|_{p,\alpha}.
 \Ee
\ecor

\subsection{Proof of Theorem \ref{thm:main6}}
For $\beta\in \{0,1/3\}$, we consider the following positive operators.
\begin{equation}\label{eq:discretoper}
Q_{\alpha,\gamma}^\beta f:=\sum_{I\in \mathcal {D}^\beta}\langle f,\frac{\chi_{Q_I}}{|I|^{2+\alpha-\gamma}}\rangle_\alpha {\chi}_{Q_I}.
\end{equation}
By comparing the positive kernel $$K_\alpha^+(z,w)=\frac{1}{|z-w|^{2+\alpha-\gamma}}$$ and the box-type kernel
$$K_\alpha^\beta(z,w):=\sum_{I\in \mathcal {D}^\beta}\frac{\chi_{Q_I}(z)\chi_{Q_I}(w)}{|I|^{2+\alpha-\gamma}},$$
one obtains the following (see \cite{PR} for the case $\gamma=0$).
\begin{proposition}\label{prop:compareP+Pbeta}
There is a constant $C>0$ such that for any $f\in L_{loc}^1{(\mathcal H)}$, $f\ge 0$, and $z\in \mathcal H$,
\begin{equation}\label{eq:compareP+gamma}
T_{\alpha,\gamma}f(z)\le C\sum_{\beta\in \{0,1/3\}}Q_{\alpha,\gamma}^\beta f(z).
\end{equation}
\end{proposition}
We can now prove Theorem \ref{thm:main6}.
\begin{proof}[Proof of Theorem \ref{thm:main6}]
It follows from the above observations that to prove the inequality (\ref{eq:main62}), it is enough to prove that under (\ref{eq:main61}), the dyadic operators $Q_{\alpha,\gamma}^\beta$ are bounded from $L^p(\mathcal{H}, \sigma^pdV_\alpha)$ into $L^q(\mathcal{H}, \omega^qdV_\alpha)$.
\vskip .3cm
We are looking to prove that there is a positive constant $C$ such that for any positive function $f$ and any $g\in L^{q'}(\mathcal{H}, dV_\alpha)$, $g\geq 0$, 
$$\int_{\mathcal{H}}(Q_{\alpha,\gamma}^\beta f(z))g(z)\omega(z)dV_\alpha(z)\le C\|\sigma f\|_{p,\alpha}\|g\|_{q',\alpha}.$$
We denote by $\overline{\Phi}$ and $\overline{\Psi}$, the complementary functions of $\Phi$ and $\Psi$ respectively. We have
\Beas
L &:=& \int_{\mathcal{H}}(Q_{\alpha,\gamma}^\beta f(z))g(z)\omega(z)dV_\alpha(z)\\ &=& \sum_{I\in \mathcal{D}^\beta}\frac{1}{|Q_I|_\alpha^{1-\frac{\gamma}{2+\alpha}}}\left(\int_{Q_I}fdV_\alpha\right)\left(\int_{Q_I}g\omega dV_\alpha\right)\\ &=& \sum_{I\in \mathcal{D}^\beta}|Q_I|_\alpha^{\frac{\gamma}{2+\alpha}}\left(\frac{1}{|Q_I|_\alpha}\int_{Q_I}(\sigma f)\sigma^{-1}dV_\alpha\right)\left(\frac{1}{|Q_I|_\alpha}\int_{Q_I}g\omega dV_\alpha\right)|Q_I|_\alpha.
\Eeas
It follows using (\ref{eq:holderhgene}) and (\ref{eq:main61}), that
\Beas
L &\le& \sum_{I\in \mathcal{D}^\beta}|Q_I|_\alpha^{\frac{\gamma}{2+\alpha}}\|\sigma f\|_{Q_I,\overline{\Phi},\alpha}\|\sigma^{-1}\|_{Q_I,\Phi,\alpha}\|g\|_{Q_I,\overline{\Psi},\alpha}\|\omega\|_{Q_I,\Psi,\alpha}|Q_I|_\alpha\\ &\le& \sum_{I\in \mathcal{D}^\beta}\|\sigma f\|_{Q_I,\overline{\Phi},\alpha}\|g\|_{Q_I,\overline{\Psi},\alpha}|Q_I|_\alpha^{\frac{1}{p}+\frac{1}{q'}}\\ &\le& \left(\sum_{I\in \mathcal{D}^\beta}\|\sigma f\|_{Q_I,\overline{\Phi},\alpha}^p|Q_I|_\alpha\right)^{1/p}\left(\sum_{I\in \mathcal{D}^\beta}\|g\|_{Q_I,\overline{\Psi},\alpha}^{p'}|Q_I|_\alpha^{\frac{p'}{q'}}\right)^{1/p'}\\ &\le& \left(\sum_{I\in \mathcal{D}^\beta}\|\sigma f\|_{Q_I,\overline{\Phi},\alpha}^p|Q_I|_\alpha\right)^{1/p}\left(\sum_{I\in \mathcal{D}^\beta}\|g\|_{Q_I,\overline{\Psi},\alpha}^{q'}|Q_I|_\alpha\right)^{1/q'}.
\Eeas
Proceeding as in the last part of the proof of Theorem \ref{thm:main5} with the help of Lemma \ref{lem:phimax}, we finally obtain
\Beas
L &:=& \int_{\mathcal{H}}(Q_{\alpha,\gamma}^\beta f(z))g(z)\omega(z)dV_\alpha(z)\\ &\le& \left(\sum_{I\in \mathcal{D}^\beta}\|\sigma f\|_{Q_I,\overline{\Phi},\alpha}^p|Q_I|_\alpha\right)^{1/p}\left(\sum_{I\in \mathcal{D}^\beta}\|g\|_{Q_I,\overline{\Psi},\alpha}^{q'}|Q_I|_\alpha\right)^{1/q'}\\ &\le& C\|\sigma f\|_{p,\alpha}\|g\|_{q',\alpha}.
\Eeas
The proof is complete.
\end{proof}
As a corollary, we have the following particular case.
\bcor\label{cor:main6}
Let $\alpha>-1$, $0\le \gamma<2+\alpha$, and $1<p\le q<\infty$.  Assume that $\omega$ is a weight and $\mu$ a positive Borel measure on $\mathcal H$ such that for some $r>1$, there is positive constant $C$ for which for any interval $I\subset \mathbb{R}$
\Be\label{eq:main61cor}
|Q_I|_\alpha^{\frac{\gamma}{2+\alpha}-\frac{1}{p}+\frac{1}{q}}\left(\frac{1}{|Q_I|_\alpha}\int_{Q_I}\omega^{-p'r}dV_\alpha\right)^{\frac{1}{p'r}}\left(\frac{1}{|Q_I|_\alpha}\int_{Q_I}\omega^{qr}dV_\alpha\right)^{\frac{1}{qr}}\le C.
\Ee

 Then there is a positive constant $K$ such that for any $f\in L^p(\mathcal H, \omega dV_\alpha)$, 
 \Be\label{eq:main62cor}
 \left(\int_{\mathcal H}\left(\omega(z)T_{\alpha,\gamma}f(z)\right)^qdV_\alpha(z)\right)^{1/q}\le K\|f\omega\|_{p,\alpha}.
 \Ee
\ecor
\subsection{Proof of Theorem \ref{thm:Apqestim}}
We start by introducing the following logarithmic maximal function (this is inspired from the definition in \cite{HyPerez}): $$\mathcal{M}_{\alpha}^{exp}f:=\sup_{I\subset \mathbb{R}}\exp\left(\frac{1}{|Q_I|_\alpha}\int_{Q_I}\log|f|dV_\alpha\right)\chi_{Q_I}.$$
It follows easily from the Jensen's inequality that $$\mathcal{M}_{\alpha}^{exp}f\le \mathcal{M}_{\alpha}f$$
consequently, $\mathcal{M}_{\alpha}^{exp}$ is bounded on $L^p(\mathcal{H}, dV_\alpha)$ for all $1<p\le \infty$. This also holds for small exponents.
\blem\label{lem:maxexpo}
Let $0<p<\infty$ and $\alpha>-1$. Then there is a positive constant $C=C(p,\alpha)$ such that
\Be\label{eq:maxexpo}
\|\mathcal{M}_{\alpha}^{exp}f\|_{p,\alpha}\le C^{1/p}\|\mathcal{M}_{\alpha}^{exp}f\|_{p,\alpha}.
\Ee
\elem
\begin{proof}
Let $0<p<q<\infty$. It is easy to see that $$\mathcal{M}_{\alpha}^{exp}f= (\mathcal{M}_{\alpha}^{exp}|f|^{p/q})^{q/p}.$$
It follows from the observations made above and the boundedness of the Hardy-Littlewood maximal in Proposition \ref{prop:sharpmaxfrac} that
\Beas
\|\mathcal{M}_{\alpha}^{exp}f\|_{p,\alpha} &:=& \|(\mathcal{M}_{\alpha}^{exp}|f|^{p/q})^{q/p}\|_{p,\alpha}\\ &=& \|(\mathcal{M}_{\alpha}^{exp}|f|^{p/q})\|_{q,\alpha}^{q/p}\\ &\le& \|(\mathcal{M}_{\alpha}|f|^{p/q})\|_{q,\alpha}^{q/p}\\ &\le& C^{1/p}\|f\|_{p,\alpha}.
\Eeas
\end{proof}
We now prove Theorem \ref{thm:Apqestim}.
\begin{proof}[Proof of Theorem \ref{thm:Apqestim}]
We note that inequality (\ref{eq:ineqApq}) is already given by Corollary \ref{cor:main44}. We then only have to prove (\ref{eq:ineqCpq}) and (\ref{eq:ineqSpq}). Once more, it is enough to check the inequality for the corresponding dyadic maximal function. Let $a\ge 2^{2+\alpha-\gamma}$. Let associate to each integer $k$, we the set
$$\Omega_{k}:=\{z\in \mathcal H: a^k<\mathcal{M}_{\alpha,\gamma}^{d,\beta}(\sigma f)(z)\leq a^{k+1}\}.$$
We already know that
$\Omega_{k}\subset \cup_{j=1}^{\infty}Q_{I_{j}}^k,$ where
$Q_{I_{j}}^k$ ($I_{j}\in \mathcal{D}^{\beta}$) is a dyadic cube maximal (with respect to the inclusion) such that
$$\frac{1}{|Q_{I_{j}}^k|_\alpha^{1-\frac{\gamma}{2+\alpha}}}\int_{Q_{I_{j}}^k}|(\sigma f)(z)|dV_\alpha(z)>a^k.$$
We start with the estimate (\ref{eq:ineqCpq}): as in the proof of the inequality (\ref{eq:main13}), we obtain
\Beas
L &:=& \int_{\mathcal H}(\mathcal{M}_{\alpha,\gamma}^{d,\beta}f(z))^{q}\omega(z)d(z)\\ &=& \sum_{k}\int_{\Omega_k}(\mathcal{M}_{\alpha,\gamma}^{d,\beta} f(z))^{q}\omega(z)d(z)\\
 &\le&
 a^{q}\sum_{k,j}\left(\frac{1}{|Q_{I_{j}}^k|_\alpha^{1-\frac{\gamma}{2+\alpha}}}\int_{Q_{I_{j}}^k}|f(z)|dV_\alpha(z)\right)^{q}|Q_{I_{j}}^k|_{\omega,\alpha}\\
 &\lesssim&
[\sigma,\omega]_{\mathcal{C}_{p,q,\alpha}}^q \sum_{k,j}\left(\int_{Q_{I_{j}}^k}|f(z)|dV_\alpha(z)\right)^{q}|Q_{I_j}^k|_{u,\alpha}^{-q/p'}\\
&\lesssim& [\sigma,\omega]_{\mathcal{C}_{p,q,\alpha}}^q\sum_{k,j}\left(\frac{1}{|Q_{I_{j}}^k|_{u,\alpha}}\int_{Q_{I_{j}}^k}|(u^{-1}f)(z)|u(z)dV_\alpha(z)\right)^{q}|Q_{I_j}^k|_{u,\alpha}^{q/p}.
\Eeas
Following again the reasoning at the end of the proof of Lemma \ref{lem:main1}, we obtain that
\Beas
\int_{\mathcal H}(\mathcal{M}_{\alpha,\gamma}^{d,\beta}f(z))^{q}\omega(z)d(z) &\lesssim& [\sigma,\omega]_{\mathcal{C}_{p,q,\alpha}}^q\left([u]_{\mathcal{B}_{\infty,\alpha}}\|u^{-1}f\|_{p,u,\alpha}^p\right)^{q/p}\\ &=& [\sigma,\omega]_{\mathcal{C}_{p,q,\alpha}}^q\left([u]_{\mathcal{B}_{\infty,\alpha}}\|\sigma f\|_{p,\alpha}^p\right)^{q/p}.
\Eeas
This completes the proof of the estimate (\ref{eq:ineqCpq}).
\vskip .3cm
Let us now prove (\ref{eq:ineqSpq}). Following the same reasoning as above , we obtain
 \Beas
L &:=& \int_{\mathcal H}(\mathcal{M}_{\alpha,\gamma}^{d,\beta}(\sigma f)(z))^{q}\omega(z)d(z)\\ 
 &\le&
 a^{q}\sum_{k,j}\left(\frac{1}{|Q_{I_{j}}^k|_\alpha^{1-\frac{\gamma}{2+\alpha}}}\int_{Q_{I_{j}}^k}|(\sigma f)(z)|dV_\alpha(z)\right)^{q}|Q_{I_{j}}^k|_{\omega,\alpha}\\
 &\lesssim&
 \sum_{k,j}\left(\frac{1}{|Q_{I_{j}}^k|_{\sigma,\alpha}}\int_{Q_{I_{j}}^k}|(\sigma f)(z)|dV_\alpha(z)\right)^{q}\frac{|Q_{I_j}^k|_{\omega,\alpha}|Q_{I_j}^k|_{\sigma,\alpha}^{q}}{|Q_{I_j}^k|_\alpha^{q(1-\frac{\gamma}{2+\alpha})}}\\
&\lesssim& [\sigma,\omega]_{\mathcal{S}_{p,q,\alpha}}\sum_{k,j}\left(\frac{1}{|Q_{I_{j}}^k|_{\sigma,\alpha}}\int_{Q_{I_{j}}^k}|(\sigma f)(z)|dV_\alpha(z)\right)^{q}|Q_{I_j}^k|_{\alpha}^{q/p}\\ && \exp\left(\frac{q}{p}\frac{1}{|Q_{I_{j}}^k|_{\alpha}}\int_{Q_{I_{j}}^k}(\log\sigma) dV_\alpha(z)\right)
\\ &\le& [\sigma,\omega]_{\mathcal{S}_{p,q,\alpha}}^{q/p}\left(\sum_{k,j}\left(\frac{1}{|Q_{I_{j}}^k|_{\sigma,\alpha}}\int_{Q_{I_{j}}^k}|(\sigma f)(z)|dV_\alpha(z)\right)^{p}\lambda_{Q_{I_j}^k}\right)^{q/p}.
\Eeas
where the sequence $\{\lambda_{Q_I}\}_{I\in \mathcal{D}^\beta}$ is defined by
$$
\lambda_{Q_I}:=\left\{ \begin{matrix} |Q_I|_{\alpha}\exp\left(\frac{1}{|Q_I|_{\alpha}}\int_{Q_{I}}(\log\sigma) dV_\alpha(z)\right) &\text{if }& Q_I=Q_{I_j}^k\,\,\,\textrm{for some}\,\,\,(k,j),\\
      0 & \text{ otherwise} .
                                  \end{matrix} \right.
$$
It follows from the Carleson embedding Theorem that 
\Beas
L &:=& \int_{\mathcal H}(\mathcal{M}_{\alpha,\gamma}^{d,\beta}(\sigma f)(z))^{q}\omega(z)d(z)\\ &\le& C\left(\int_{\mathcal{H}}|f(z)|^p\sigma(z)dV_\alpha(z)\right)^{q/p}
\Eeas
provide $\{\lambda_{Q_I}\}_{I\in \mathcal{D}^\beta}$ is a $(\sigma,\alpha,1)$-Carleson sequence.  Let us check the latter. For any interval $J\in \mathcal{D}^\beta$, we have using Lemma \ref{lem:maxexpo},
\Beas 
S &:=& \sum_{I\subseteq J, I\in \mathcal{D}^\beta}|Q_I|_{\alpha}\exp\left(\frac{1}{|Q_I|_{\alpha}}\int_{Q_{I}}(\log\sigma) dV_\alpha(z)\right)\\ &\simeq& \sum_{I\subseteq J, I\in \mathcal{D}^\beta}|T_I|_{\alpha}\exp\left(\frac{1}{|Q_I|_{\alpha}}\int_{Q_{I}}(\log\sigma) dV_\alpha(z)\right)\\ &=& \sum_{I\subseteq J, I\in \mathcal{D}^\beta}\int_{T_I}\exp\left(\frac{1}{|Q_I|_{\alpha}}\int_{Q_{I}}(\log\sigma) dV_\alpha(z)\right)dV_\alpha\\ &\le& \sum_{I\subseteq J, I\in \mathcal{D}^\beta}\int_{T_I}\mathcal{M}_\alpha^{exp}(\sigma \chi_{Q_J})dV_\alpha\\ &\le& \int_{Q_J}\mathcal{M}_\alpha^{exp}(\sigma \chi_{Q_J})dV_\alpha\\ &\le& C_\alpha |Q_J|_{\sigma,\alpha}.
\Eeas
The proof is complete.

\end{proof}
\subsection{Proof of Theorem \ref{thm:Cpqestim}}
For the proof of the inequality (\ref{eq:ineqCpqsharp1}), it is enough to prove the following.
\bprop
Let $\alpha>-1$,  $0\le \gamma<2+\alpha$, and $1<p\le q<\infty$.  Let $\sigma,\omega$ be weights and put $u=\sigma^{-p'}$. Then

\Be\label{eq:ineqCpqsharpdyadic}
\|\mathcal{M}_{\alpha,\gamma}^{d,\beta} f\|_{q,\omega,\alpha}\lesssim [\sigma,\omega]_{\mathcal{C}_{p,q,\alpha}}[u]_{\mathcal{B}_{\infty,\alpha}}^{\frac{1}{q}}\|\sigma f\|_{p,\alpha}.
\Ee
\eprop
\begin{proof}
We use the same notations as in the proof of Theorem \ref{thm:Apqestim}. Let $\theta$
be such that $\frac{1}{p}-\frac{1}{q}=\frac{\theta}{2+\alpha}$. Then with the squares $Q_{I_{j}}^k$ as above, we obtain
\Beas
L &:=& \int_{\mathcal H}(\mathcal{M}_{\alpha,\gamma}^{d,\beta}f(z))^{q}\omega(z)d(z)\\ 
 &\lesssim&
 \sum_{k,j}\left(\frac{1}{|Q_{I_{j}}^k|_\alpha^{1-\frac{\gamma}{2+\alpha}}}\int_{Q_{I_{j}}^k}|f(z)|dV_\alpha(z)\right)^{q}|Q_{I_{j}}^k|_{\omega,\alpha}\\
 &\lesssim&
[\sigma,\omega]_{\mathcal{C}_{p,q,\alpha}}^q \sum_{k,j}\left(\int_{Q_{I_{j}}^k}|f(z)|dV_\alpha(z)\right)^{q}|Q_{I_j}^k|_{u,\alpha}^{-q/p'}\\
&\lesssim& [\sigma,\omega]_{\mathcal{C}_{p,q,\alpha}}^q\sum_{k,j}\left(\frac{1}{|Q_{I_{j}}^k|_{u,\alpha}^{1-\frac{\theta}{2+\alpha}}}\int_{Q_{I_{j}}^k}|(u^{-1}f)(z)|u(z)dV_\alpha(z)\right)^{q}|Q_{I_j}^k|_{u,\alpha}.
\Eeas
As $\{|Q_{I}|_{u,\alpha}\}_{I\in \mathcal{D}^\beta}$ is a $(u,\alpha,1)$-Carleson sequence with Carleson constant $[u]_{\mathcal{B}_{\infty,\alpha}}$, we obtain using Theorem \ref{thm:Carlembed0} and Corollary \ref{cor:weightedfracmax} that 
\Beas
L_1 &:=& \sum_{k,j}\left(\frac{1}{|Q_{I_{j}}^k|_{u,\alpha}^{1-\frac{\theta}{2+\alpha}}}\int_{Q_{I_{j}}^k}|(u^{-1}f)(z)|u(z)dV_\alpha(z)\right)^{q}|Q_{I_j}^k|_{u,\alpha}\\ &\le& [u]_{\mathcal{B}_{\infty,\alpha}}\|\mathcal{M}_{u,\alpha,\theta}^{d,\beta}(u^{-1}f)\|_{q,u,\alpha}^q\\ &\lesssim& [u]_{\mathcal{B}_{\infty,\alpha}}\|u^{-1}f\|_{p,u,\alpha}^q= [u]_{\mathcal{B}_{\infty,\alpha}}\|\sigma f\|_{p,\alpha}^q.
\Eeas
Hence
$$\int_{\mathcal H}(\mathcal{M}_{\alpha,\gamma}^{d,\beta}f(z))^{q}\omega(z)d(z)\lesssim [\sigma,\omega]_{\mathcal{C}_{p,q,\alpha}}^q[u]_{\mathcal{B}_{\infty,\alpha}}\|\sigma f\|_{p,\alpha}^q.$$
The proof is complete.

\end{proof}
\subsection{Proof of Theorem \ref{thm:sharpmaxfrac}}
We observe again that to prove the estimate (\ref{eq:ineqBpq}), it is enough to prove the following.
\bprop\label{prop:sharpmaxfrac}
Let $\alpha>-1$,  $0\le \gamma<2+\alpha$, and $1<p<\frac{2+\alpha}{\gamma}$. Suppose $q$ is defined by the relation $\frac{1}{q}=\frac{1}{p}-\frac{\gamma}{2+\alpha}$. If $\omega\in \mathcal{B}_{p,q,\alpha}$, then 
\Be\label{eq:dyadicineqBpq}
\|\omega \mathcal{M}_{\alpha,\gamma}^{d,\beta} f\|_{q,\alpha}\le ([\omega]_{\mathcal{B}_{p,q,\alpha}})^{\frac{p'}{q}(1-\frac{\gamma}{2+\alpha})}\|\omega f\|_{p,\alpha}.
\Ee
\eprop
\begin{proof}[Proof of Proposition \ref{prop:sharpmaxfrac}]
We use an idea from \cite{Laceyetal}. We recall that for a weight $\sigma$, $\mathcal{M}_{\sigma,\alpha,\gamma}$ is the weighted fractional maximal function as defined in (\ref{eq:defweightedfracmax}). When $\gamma=0$, we write $\mathcal{M}_{\sigma,\alpha}$ for the corresponding weighted Hardy-Littlewood maximal function. Let us put $u=\omega^q$, and $v=\omega^{-p'}$. Define $r=1+\frac{q}{p'}$ and observe that its conjugate exponent is $r'=1+\frac{p'}{q}=p'(1-\frac{\gamma}{2+\alpha})$. For any dyadic interval $I\subset \mathbb{R}$, we first obtain
\Beas
S &:=&\frac{1}{|Q_I|_\alpha^{1-\frac{\gamma}{2+\alpha}}}\int_{Q_I}|f|dV_\alpha\\ &=& \left(\frac{|Q_I|_{u,\alpha}|Q_I|_{v,\alpha}^{q/p'}}{|Q_I|_\alpha^{\frac{q}{p'}+1}}\right)^{\frac{r'}{q}}\left(\frac{|Q_I|_\alpha}{|Q_I|_{u,\alpha}}\right)^{\frac{r'}{q}}\left(\frac{1}{|Q_I|_{v,\alpha}^{1-\frac{\gamma}{2+\alpha}}}\int_{Q_I}|f|dV_\alpha\right)\\ &\le& ([\omega]_{\mathcal{B}_{p,q,\alpha}})^{\frac{r'}{q}}\left(\frac{|Q_I|_\alpha}{|Q_I|_{u,\alpha}}\right)^{\frac{r'}{q}}\left(\frac{1}{|Q_I|_{v,\alpha}^{1-\frac{\gamma}{2+\alpha}}}\int_{Q_I}|v^{-1}f|vdV_\alpha\right)\\ &\le& ([\omega]_{\mathcal{B}_{p,q,\alpha}})^{\frac{r'}{q}}\left(\frac{|Q_I|_\alpha}{|Q_I|_{u,\alpha}}\right)^{\frac{r'}{q}}\mathcal{M}_{v,\alpha,\gamma}^d(v^{-1}f)\\ &=& ([\omega]_{\mathcal{B}_{p,q,\alpha}})^{\frac{r'}{q}}\left(\frac{|Q_I|_\alpha}{|Q_I|_{u,\alpha}}(\mathcal{M}_{v,\alpha,\gamma}^d(v^{-1}f))^{\frac{q}{r'}}\right)^{\frac{r'}{q}}\\ &\le& ([\omega]_{\mathcal{B}_{p,q,\alpha}})^{\frac{r'}{q}}\left(\frac{1}{|Q_I|_{u,\alpha}}\int_{Q_I}(\mathcal{M}_{v,\alpha,\gamma}^d(v^{-1}f))^{\frac{q}{r'}}dV_\alpha\right)^{\frac{r'}{q}}\\ &\le& ([\omega]_{\mathcal{B}_{p,q,\alpha}})^{\frac{r'}{q}}\left(\mathcal{M}_{u,\alpha}^d\left((\mathcal{M}_{v,\alpha,\gamma}^d(v^{-1}f))^{\frac{q}{r'}}u^{-1}\right)\right)^{\frac{r'}{q}}
\Eeas
Now, using Corollary \ref{cor:weightedfracmax}, we obtain
\Beas
\|\omega \mathcal{M}_{\alpha,\gamma}^d f\|_{q,\alpha} &=& \|\mathcal{M}_{\alpha,\gamma}^d f\|_{q,u,\alpha}\\ &\le& ([\omega]_{\mathcal{B}_{p,q,\alpha}})^{\frac{r'}{q}}\|\mathcal{M}_{u,\alpha}^d\left((\mathcal{M}_{v,\alpha,\gamma}^d(v^{-1}f))^{\frac{q}{r'}}u^{-1}\right)\|_{r',u,\alpha}^{r'/q}\\ &\le& C([\omega]_{\mathcal{B}_{p,q,\alpha}})^{\frac{r'}{q}}\|(\mathcal{M}_{v,\alpha,\gamma}^d(v^{-1}f))^{\frac{q}{r'}}u^{-1}\|_{r',u,\alpha}^{r'/q}\\ &=& C([\omega]_{\mathcal{B}_{p,q,\alpha}})^{\frac{r'}{q}}\|\mathcal{M}_{v,\alpha,\gamma}^d(v^{-1}f)\|_{q,v,\alpha}\\ &\le& C([\omega]_{\mathcal{B}_{p,q,\alpha}})^{\frac{r'}{q}}\|v^{-1}f\|_{p,v,\alpha}= C([\omega]_{\mathcal{B}_{p,q,\alpha}})^{\frac{r'}{q}}\|\omega f\|_{p,\alpha}.
\Eeas
The proof is complete.
\end{proof}
\section{Example}
We give examples to show that the constants in (\ref{eq:ineqCpqsharp}) and (\ref{eq:ineqBpq}) are sharp.
We start with (\ref{eq:ineqBpq}): we recall that $0\le \gamma <2+\alpha$. Fix $0<\epsilon<1$.  We consider $\omega(z)=|z|^{\frac{(2+\alpha-\epsilon)}{p'}}$. One easily check that $\omega\in \mathcal{B}_{p,q,\alpha}$ and that $[\omega]_{\mathcal{B}_{p,q,\alpha}}\simeq \epsilon^{-\frac{q}{p'}}$. We also consider the function $f(z)=|z|^{\epsilon-(2+\alpha)}{\bf 1}_{\{z\in \mathcal{H}:|z|\leq 1\}}$. We obtain that $\|\omega f\|_{p,\alpha}\simeq \epsilon^{-1/p}$.
\vskip .3cm
Let $z\in \{z\in \mathcal{H}:|z|\leq 1\}$. Then
\Beas
\mathcal{M}_{\alpha,\gamma}f(z) &\geq& \frac{1}{|z|^{2+\alpha-\gamma}}\int_{\{w\in \mathcal{H}:|w|\leq |z|\}}|f(w)|dV_\alpha(w)\\ &\geq& \frac{1}{|z|^{2+\alpha-\gamma}}\int_{\{w\in \mathcal{H}:|w|\leq |z|\}}|w|^{\epsilon-2-\alpha}dV_\alpha(w)\\ &\simeq& |z|^{\epsilon-2-\alpha+\gamma}\epsilon^{-1}.
\Eeas
Hence
\Beas
\int_{\mathcal{H}}(\omega(z)\mathcal{M}_{\alpha,\gamma}f(z))^qdV_\alpha(z) &\geq& C_{\alpha,q}\epsilon^{-q}\int_{\{z\in \mathcal{H}:|z|\leq 1\}}|z|^{q(\epsilon-2-\alpha+\gamma)+\frac{q}{p'}(2+\alpha-\epsilon)}dV_\alpha(z)\\ &=& C_{\alpha,q}\epsilon^{-q}\int_{\{z\in \mathcal{H}:|z|\leq 1\}}|z|^{-2-\alpha+\frac{q}{p}\epsilon}dV_\alpha(z)\\ &\simeq& \epsilon^{-q-1}.  
\Eeas
Thus 
\Beas 
\|\omega \mathcal{M}_{\alpha,\gamma}f\|_{q,\alpha} &\gtrsim& \epsilon^{-(\frac{1}{p'}+\frac{1}{q})}\epsilon^{-\frac{1}{p}}\simeq [\omega]_{\mathcal{B}_{p,q,\alpha}}^{(1-\frac{\gamma}{2+\alpha})\frac{p'}{q}}\|\omega f\|_{p,\alpha},
\Eeas
showing that (\ref{eq:ineqBpq}) is sharp.
\vskip .3cm
Let us also check that (\ref{eq:ineqCpqsharp}) is sharp. Fix $0<\epsilon<1$.  We consider the same weight and function as above, $\omega(z)=|z|^{\frac{(2+\alpha-\epsilon)}{p'}}$ and $f(z)=|z|^{\epsilon-(2+\alpha)}{\bf 1}_{\{z\in \mathcal{H}:|z|\leq 1\}}$. Recall that $u=\omega^{-p'}$. We obtain $[u]_{\mathcal{B}_{\infty,\alpha}}\le [u]_{\mathcal{B}_{p',\alpha}}\simeq \frac{1}{\epsilon}$. From the previous computations, we have
\Beas 
\|\omega \mathcal{M}_{\alpha,\gamma} f\|_{q,\alpha} &\gtrsim& \epsilon^{-1-\frac{1}{q}}\\ &=&\epsilon^{-\frac{1}{p'}}\epsilon^{-\frac{1}{q}}\epsilon^{-\frac{1}{p}}\\ &\gtrsim& [\omega]_{\mathcal{B}_{p,q,\alpha}}^{\frac{1}{q}}[u]_{\mathcal{B}_{\infty,\alpha}}^{\frac{1}{q}}\|\omega f\|_{p,\alpha},
\Eeas
proving the sharpness of (\ref{eq:ineqCpqsharp}).

\end{document}